\documentclass[12pt,reqno]{article}

\title{Growth and Percolation on the Uniform Infinite Planar Triangulation}
\author{Omer Angel}
\date{August 15, 2002}

\usepackage{latexsym, amsfonts, amsmath, amsthm, pstricks}

\numberwithin{equation}{section}
\numberwithin{figure}{section}

\newtheorem{thm}{Theorem} [section]
\newtheorem{prop}[thm]{Proposition}
\newtheorem{coro}[thm]{Corollary}
\newtheorem{conj}[thm]{Conjecture}
\newtheorem{lemma}[thm]{Lemma}
\newtheorem{claim}[thm]{Claim}
\newtheorem{defn}[thm]{Definition}
\newtheorem{fact}[thm]{Fact}

\newcommand{\E}{{\mathbb E}}
\newcommand{\F}{{\mathcal F}}
\newcommand{\N}{{\mathbb N}}
\newcommand{\R}{{\mathbb R}}
\newcommand{\barB}{{\overline{B}}}
\newcommand{\barX}{{\overline{X}}}
\renewcommand{\P}{{\mathbb P}}
\newcommand{\D}{\partial}
\renewcommand{\phi}{\varphi}
\newcommand{\ep}{{\varepsilon}}
\newcommand{\slt}{\prec}

\newcommand{\note}{{\noindent \bf Note.\ }}
\newcommand{\pexp}[1]{\exp \left( #1 \right)}

\newcommand{\be}{\begin{enumerate}}
\newcommand{\ee}{\end{enumerate}}
\newcommand{\bi}{\begin{itemize}}
\newcommand{\ei}{\end{itemize}}

\begin{document}

\maketitle

\begin{abstract}
A construction as a growth process for sampling of the uniform infinite
planar triangulation (UIPT), defined in \cite{UIPT1}, is given. The
construction is algorithmic in nature, and is an efficient method of
sampling a portion of the UIPT. 

By analyzing the progress rate of the growth process we show that a.s.\ the
UIPT has growth rate $r^4$ up to polylogarithmic factors, confirming
heuristic results from the physics literature. Additionally, the boundary
component of the ball of radius $r$ separating it from infinity a.s.\ has
growth rate $r^2$ up to polylogarithmic factors. It is also shown that the
properly scaled size of a variant of the free triangulation of an $m$-gon
(also defined in \cite{UIPT1}) converges in distribution to an asymmetric
stable random variable of type $1/2$.

By combining Bernoulli site percolation with the growth process for the
UIPT, it is shown that a.s.\ the critical probability $p_c=1/2$ and that at
$p_c$ percolation does not occur.
\end{abstract}


%%%%%%%%%%%%%%%%%%%%%%
\section{Introduction}  \label{sec:intro}
%%%%%%%%%%%%%%%%%%%%%%

Since the 1960's there has been a combinatorial study of the properties of
finite planar maps chosen uniformly among members of one of a number of
classes. Many of their statistical properties are known. Tutte \cite{Tut2}
has shown that almost all planar map have no non-trivial symmetries. Later
by other means this result was extended to many other classes of planar
structures \cite{RiWo2}. 

The classes of planar objects contain the class of general planar maps as
well as the classes of triangulations (where all faces are triangles),
quadrangulations and maps with other possible faces, and others. Members of
a class may also have restrictions on degrees of vertices or on their
connectivity. While understanding each particular class involves different
techniques for handling the particular difficulties the class presents,
many results appear to hold universally for any class of planar structures
examined.

Another aspect common to many of the statistical results regarding
uniformly chosen planar structures is dealing with asymptotics of the
structure's properties as the size tends to infinity. Thus there are
results about the distribution of degrees in a uniformly chosen
triangulation \cite{GaRi}, the size of 3-connected components \cite{BRW},
and probabilistic 0-1 laws \cite{BCR}. Many of these results may be viewed
as a finite version of corresponding results regarding an infinite planar
structure of the same class. Considering infinite versions of such results
may make the results more concise, though some of the precision of the
finite case may be lost. For example, asymmetry of all but an exponentially
small portion of planar maps translates to the simple (but no longer
quantitative) fact that the infinite object a.s.\ has no symmetries.

\medskip

The uniform infinite planar triangulation (UIPT) is one case of such an
infinite object. It is defined by considering the uniform measure on planar
triangulations with $n$ vertices and taking a weak limit as $n$ tends to
infinity. This gives rise to a probability measure supported on infinite
planar triangulations, with the UIPT denoting a sample. The UIPT was
suggested in \cite{BeSc} and shown to exist in \cite{UIPT1} where some of
its properties are studied. Of course a similar limit may be taken for any
class of planar structures.

We consider triangulations were chosen solely because they appear to be
easier to work with then other classes, and yield better to the techniques
in this paper. However, this is not a crucial necessity. While local
properties such as degree distribution are dependent on the class, large
scale properties appear to be independent of local definitions. This is
demonstrated in the infinite setting in \cite{UIPT1} where two types of
triangulations are studied --- with or without double edges --- and a
simple relation between them is given. This relation implies that on a
large scale the two objects have the same properties.

\medskip

Schaeffer~\cite{Scha} found a bijection between certain types planar maps
and labeled trees. Chassaing and Schaeffer~\cite{ChSc} recently used this
bijection to show a connection between the asymptotic distribution of the
radius of a random map and the integrated super-Brownian excursion. They
deduce from this connection that the diameter of such a map of size $n$ 
scales as $n^{1/4}$. The results presented here on the volume growth of the
UIPT are a kind of infinite version of their result.

While they work with planar quadrangulations and we 
with triangulations, it appears that such local differences are
insignificant when large scale observations such as diameter, growth,
separation, etc. are concerned. This phenomena is referred to as
universality.

Physicists study similar random planar structures under the title of {\em
2-dimensional quantum gravity}. The essential idea here is to develop a
quantum theory of gravity by extending to higher dimensions the concept of
Feynman integrals for one dimensional paths. Triangulations (or other
classes of maps) are viewed as a discretized version of a 2 dimensional
manifold. In two dimensions this gives rise to a rich theory, much of which
appears to be missing for higher dimensions. More often physicists are
interested not in the discretized planar triangulation but in a continuous
scaling limit of it which is believed to exist.

Once the basic structure is defined, models of statistical physics may be
introduced in it and in many cases solved (e.g. \cite{BoKa} and others).
Physicists applied here the methods of random matrix models \cite{FGZJ}.
Through these methods and other heuristics many conjectures were made on
the structure of such triangulations. In particular, it is believed that
the Hausdorff dimension of the scaling limit of 2-dimensional quantum
gravity is 4 \cite{AmWa}. This is a continuous form of the volume growth
results of Sections~\ref{sec:hull_vol} and \ref{sec:ball_vol}. For a good
general exposition of quantum gravity see \cite{ADJ}, as well as
\cite{Amb,Dav}.

Another recent version of the growth results may be found in \cite{ChSc}.
Schaeffer~\cite{Scha} constructed a bijection between certain types planar
maps and labeled trees \cite{Scha}. Chassaing and Schaeffer were able to
use that bijection to establish a connection between the asymptotic
distribution of the radius of a random map and the integrated
super-Brownian excursion. They deduce from this connection that the
diameter of such a map of size $n$ scales as $n^{1/4}$.

\medskip

A further important motivation for understanding random triangulations (and
random planar structures in general) stems from the KPZ relation
\cite{KPZ}. On a random surface many models of statistical physics become
easier to analyze then in the Euclidean plane, as some of the geometric
aspects of the problem disappear or can be disregarded. The KPZ relation,
while not rigorously understood, is a relation between critical exponents
of models on a random planar surface and the corresponding exponents in the
plane.

For example non-intersection exponents for Brownian motion in the
Euclidean plane or half plane correspond to asymptotic non-intersection
exponents for random walks on a random surface. Similar relations hold for
exponents governing behavior of self avoiding or loop erased random walks,
boundary geometry of clusters in percolation, Ising or Potts model, and
more. Using this relation the values of Brownian motion intersection
exponents were calculated \cite{KPZ1,KPZ2}. Later a rigorous derivation of
the same values was found using the $SLE$ process \cite{LSW1,LSW2,LSW3}.

\medskip

The general structure of this paper is as follows:
In Section~\ref{sec:peeling} a process for sampling the UIPT is presented
--- the peeling process. A Representation of the UIPT as the product of a
growth process with a relatively simple definition and steps makes possible
the analysis in the following sections.

In Section~\ref{sec:markov} a key aspect of the sampling process is
considered --- the boundary size of the finite triangulation generated
after a number of steps. This boundary size, apart from being independently
interesting in connection with separation, proves to be essential for
understanding the growth of the UIPT and (to a lesser extent) Bernoulli
percolation.

In Sections~\ref{sec:boundary}, \ref{sec:hull_vol}, \ref{sec:ball_vol}
respectively asymptotic results on the boundary size of the ball, the
hull's volume and the ball's volume are proved.

Finally, in Section~\ref{sec:percolation}, percolation of the UIPT is
studied. The analysis here is based on a significant simplification of 
Bernoulli percolation derived from the construction of
Section~\ref{sec:peeling}, and depends only weekly on results in other
parts of the paper.

We now proceed to state the main results of the paper, followed by some
needed background and results (both general and specific).

\medskip

A note on notation: By $X_n \sim Y_n$ we mean that $\frac{X_n}{Y_n} \to 1$.
By $X_n\approx Y_n$ we mean that $\frac{\log X_n}{\log Y_n} \to 1$.

\medskip

The author thanks Oded Schramm, Itai Benjamini and B\`alint Vir\`ag for
helpful discussions. Part of this research was made during visits of the
author to Microsoft research. The author thanks his host for those visits.

%%%%%%%%%%%%%%%%%%%%%%%%%
\subsection{Main Results}
%%%%%%%%%%%%%%%%%%%%%%%%%

We consider in this paper only type II triangulations, i.e.\ planar
triangulations with possibly double edges but no loops. The results on
growth and percolation may be translated to type III triangulations through
the relation between the two UIPT laws \cite{UIPT1}, and to type I
triangulations through a similar decomposition. Precise definitions of
these classes appear in \cite{UIPT1}.

The UIPT is the law of a measure on infinite rooted planar graphs. Let
$B_r$ be the ball of radius $r$ (w.r.t.\ the graph metric) around the root
in the UIPT, thus $B_r$ is a finite sub-triangulation of the UIPT. The
complement of the ball is not generally connected. Denote by $\barB_r$ the
hull of the ball consisting of $B_r$ together with all finite components of
the complement. If $|T|$ is the number of vertices in a triangulation $T$,
then we prove

\begin{thm} \label{thm:hull_growth}
For any $\ep>0$, a.s.
\[
\limsup_{r \to \infty} \frac{|\barB_r|}{r^4 \log^{6+\ep} r} < \infty ,
\]
and a.s.
\[
\lim_{r \to \infty} \frac{|\barB_r| \log^{32/3+\ep} r}{r^4} = \infty .
\]
\end{thm}

We also prove for the ball itself:

\begin{thm} \label{thm:ball_growth}
A.s.
\[
\limsup_{r \to \infty} \frac{|B_r|}{r^4 \log^{6+\ep} r} < \infty ,
\]
and for any $\ep>0$
\[
\lim_{r \to \infty} \frac{|B_r| \log^{33/2+\ep} r}{r^4} = \infty .
\]
\end{thm}

The identical upper bounds appear since one is a corollary of the other.
While for the lower bounds the reduction works in the opposite direction,
$|\barB_r|$ is easier to analyze as an intermediate stage giving the
slightly better lower bound. However, this is of little significance as
in the above (and following) results the powers of the logarithms are
probably not the best possible. In proving these theorems the quantities in
question are expressed as cumulative sums of random variables which are
neither independent nor identically distributed, but are independent enough
for some methods. Thus the above results are a sort of LIL for these sums.
The proof as well as other recent results \cite{ChSc} suggest

\begin{conj}
The random variables $r^{-4}|B_r|$ and $r^{-4}|\barB_r|$ converge in
distribution. 
\end{conj}

This is roughly a converse of the result of \cite{ChSc}, stating that the
radius of a uniform planar quadrangulation of size $N$ scales as $N^{1/4}$.

While it is the ball growth rate of graphs which historically was the focus
of more extensive research, the hull is also of independent interest (even
apart from its submission to analysis). For one thing, there are questions
regarding separation properties and the isoperimetric inequality in random
planar maps. Hulls of balls (not necessarily around the root vertex) are
candidates for having small boundary sizes compared to their volume.

Let $\D \barB_r$ denote the outer boundary of $B_r$, i.e. the vertices with
neighbors in the infinite component of the complement. For $\D \barB_r$ we
prove the following:

\begin{thm} \label{thm:boundary_size}
A.s. the size of the outer boundary of the ball of radius $r$, satisfies
\[
\limsup_{r \to \infty} \frac{|\D \barB_r|}{r^2 \log^3 r} < \infty ,
\]
and for any $\ep>0$
\[
\lim_{r \to \infty} \frac{|\D \barB_r| \log^{6+\ep} r}{r^2} = \infty .
\]
\end{thm}

In particular this demonstrates that the hull's boundary is roughly the
square root of their volume. While generally a large set with a small
boundary size will not be a ball centered at the root vertex, it is
plausible that it is not very different from a ball around some vertex.
This suggests an anchored isoperimetric inequality, saying that the minimal
boundary size for a connected set $S$ in the UIPT of size $n$ containing
the root $\min_{0\in S, |S|=n} |\D S|$ scales as $\sqrt n$.

This is in contrast to the finite case where a heuristic argument suggests
a possible scaling of $n^{1/4}$. The argument is that in a uniform map of
size $n$ typical distances are on the order of $n^{1/4}$ (see \cite{ChSc}
and the argument of Section~\ref{sec:ball_vol}). Hence it is likely that to
separate such a map into two roughly equal part one needs a cycle of length
of order $n^{1/4}$. Of course, there is no need for the two exponents to
coincide, as the finite version of the anchored isoperimetric inequality
for the UIPT is to find the minimal boundary of a set of fixed size
containing the root, and not of a set of roughly half the vertices.

\medskip

When the UIPT is sampled using the growth process it is easy to add random
colors to the vertices. This results in a sample of Bernoulli site
percolation on the UIPT. Hence, percolation on the UIPT is effectively
reduced to a simple Markov chain. Using this approach we see that as far as
percolation goes the UIPT is similar to the triangular lattice in $\R^2$:

\begin{thm} \label{thm:percolation}
On the UIPT, the critical probability for site percolation is a.s.\ 1/2.
Moreover, at $p=1/2$ a.s.\ there are no infinite clusters.
\end{thm}

A part of the proof is a 0-1 result for the UIPT and for percolation on it
that is a consequence of the peeling construction.

%%%%%%%%%%%%%%%%%%%%%%%%%%%%%%%
\subsection{Further Background}
%%%%%%%%%%%%%%%%%%%%%%%%%%%%%%%

The following results either appear in \cite{UIPT1} and are repeated here
for completeness or are general facts needed.

%%%%%%%%%%%%%%%%%%%%%%%%%%%%%%%%%%%%%%%
\subsubsection{Counting Triangulations}
%%%%%%%%%%%%%%%%%%%%%%%%%%%%%%%%%%%%%%%

Many of the results derived about the UIPT are in essence consequences of
the asymptotics of the formula for the number of triangulations of a given
size. These asymptotics are common to many other classes of planar
structures. The following combinatorial result may be found in
\cite{GFBook}. It is derived using the techniques introduced by Tutte
\cite{Tutte}.

\begin{prop}\label{prop:count}
For $n,m\geq 0$, not both 0, the number of rooted type II triangulations of
a disc with $m+2$ boundary vertices and $n$ internal vertices is
\[
\phi_{n,m} = \frac {2^{n+1} (2m+1)! (2m+3n)!} {m!^2 n! (2m+2n+2)! } .
\]
\end{prop}

The case $n=m=0$ for type II triangulations warrants special attention. A
triangulation of a 2-gon must have at least one internal vertex so there
are no triangulations with $n=m=0$, yet the above formula gives
$\phi_{0,0}=1$. It will be convenient to use this value rather than 0 for 
the following reason: Typically a triangulation of an $m$-gon is used not
in itself but is used to close an external face of size $m$ of some other
triangulation by ``gluing'' it in. When the external face is a 2-gon, there 
is a further possibility of closing the hole by gluing the two edges to
each other with no additional vertices. Setting $\phi_{0,0}=1$ takes this 
possibility into account. The formula will therefore be used also for
$n=m=0$.

Using Stirling's formula, the asymptotics of this are found to be
\[
\phi_{n,m} \sim C_m \alpha^n n^{-5/2} ,
\]
where $\alpha_2=27/2$ and
\[
C_m = \frac{\sqrt{3}(2m+1)!}{4\sqrt{\pi}m!^2} (9/4)^m
                \sim C 9^m m^{1/2} .
\]
The power terms $n^{-5/2}$ and $m^{1/2}$ are common to many classes of
planar structures. They arise from the common observation that a cycle
partitions the plane into two parts (Jordan's curve Theorem) and that the
two parts may generally be triangulated (or for other classes, filled)
independently of each other.

We also define and use the partition function for triangulations of an
$(m+2)$-gon:
\[
Z_m(t) = \sum \phi_{n,m} t^n .
\]

The following appears (up to a reparametrization) in \cite{GFBook} as an
intermediate step in the derivation of Proposition~\ref{prop:count}.

\begin{prop}\label{prop:Z}
If $t=\theta(1-2\theta)^2$ with $\theta \in [0,1/6]$, then
\[
Z_m(t) = \frac{(2m)!((1-6\theta) m +2-6\theta)}{m!(m+2)!}
           (1-2\theta)^{-(2m+2)} .
\]
\end{prop}

As a corollary, by setting $\theta=1/6$ and $t=2/27=\alpha^{-1}$ we get:
\[
Z_m = Z_m(\alpha^{-1}) =
      \frac{(2m)!}{m!(m+2)!} \left(\frac{9}{4}\right)^{m+1}.
\]

%%%%%%%%%%%%%%%%%%%%%%%%
\subsubsection{Locality}
%%%%%%%%%%%%%%%%%%%%%%%%

The free triangulation of a disc is defined in \cite{UIPT1}:

\begin{defn} \label{def:free}
The {\em free distribution} on rooted triangulations of an $(m+2)$-gon,
denoted $\mu_m$, is  the probability measure that assigns weight 
\[
 \alpha^{-n} / Z_m(\alpha^{-1})
\]
to each rooted triangulation of the $(m+2)$-gon having $n$ internal
vertices.
\end{defn}

One more distribution on triangulations of a disc is given by the UIPT of a
disc. The limit of uniform distributions on triangulations of the sphere
with $N$ vertices as $N \to \infty$ exists. By conditioning on the root
having $m$ distinct neighbors with a single edge to each, and removing
the triangles incident on the root (and choosing a new root) an immediate
corollary is that there also exists a uniform distribution on infinite
triangulations of an $m$-gon.

\begin{defn} \label{def:rigid}
A rooted triangulation $A$ is {\em rigid} if no triangulation includes two
distinct copies of $A$ with coinciding roots.
\end{defn}

In particular, if the triangles of a triangulation form a connected set in
the dual graph of the triangulation, then the triangulation is rigid. This
is not a necessary condition, but it suffices for current needs.

For a finite rooted triangulation $A$ with $k$ external faces, define the
event $R_i(A)$ as the set of all infinite rooted triangulations $T$ of the
plane that include $A$ as a sub-triangulation with the roots coinciding,
and such that the component of $T$ in the $i$'th face of $A$ is
infinite. In \cite{UIPT1} it is shown that $\tau(R_i(A) \cap R_j(A))=0$ for
any $i\neq j$.

A result of \cite{UIPT1} that is a basic tool for the construction of
the next section is the following:

\begin{thm}\label{thm:locality}
Let $A$ be a finite rigid triangulation. Assume $A$ has $k$ external faces
of sizes $m_1+2,\ldots,m_k+2$. Condition on the event $R_i(A)$, and let
$T_j$ denote the component of the UIPT in the $j$'th face. Then: 
\be
\item
The triangulations $T_j$ are independent.
\item
$T_i$ has the same law as the UIPT of an $(m_i+2)$-gon.
\item
For $j \ne i$, $T_j$ has the same law as the free triangulation of an
$(m_j+2)$-gon. 
\ee
\end{thm}

%%%%%%%%%%%%%%%%%%%%%%%%%%%%%%%%%%%%%%%
\subsubsection{Stable random variables}
%%%%%%%%%%%%%%%%%%%%%%%%%%%%%%%%%%%%%%%

For any $\alpha \in (0,2)$ a completely asymmetric stable random variable
of type $\alpha$ will be denoted by $S_\alpha$. These are real random
variables with the property that the sum of $n$ i.i.d.\ copies of
$S_\alpha$ is distributed like $n^{1/\alpha} S_\alpha$. The completely
asymmetric stable random variables are characterized by having density
functions with super-polynomial decay on the left (as in the last three
properties below). We will need the following facts about stable random
variables (see \cite{SaTa,Zolo}):

\begin{fact}\label{fact:stable_dist}
\be
\item As $t \to \infty$, $\P(S_{\alpha}>t) \sim c t^{-\alpha}$.
\item For $0<\alpha<1$, a.s.\ $S_\alpha>0$.
\item For $0<\alpha<1$ there exists $c>0$ such that for small positive $t$,\\
$\P(S_\alpha<t) < \pexp{-ct^{\alpha/(\alpha-1)}}$.
\item For $1<\alpha<2$ there exists $c>0$ such that for $t>0$,\\
$\P(S_\alpha<-t) < \pexp{-ct^{\alpha/(\alpha-1)}}$.
\ee
\end{fact}

Of particular interest is $S_{3/2}$ having the Airy distribution, with
noted connections to random planar structures \cite{BFSS1,BFSS2}, and
$S_{1/2}$ having the Levi distribution.

%%%%%%%%%%%%%%%%%%%%%%%%%%
\section{A Growth Process} \label{sec:peeling}
%%%%%%%%%%%%%%%%%%%%%%%%%%

A possible method of sampling the UIPT is by adding one triangle at a
time to a finite sub-triangulation, each time adding a triangle with the
appropriate distribution conditioned on the sub-triangulations sampled so
far. This has some advantages over adding all the vertices a given distance
from the root at once. Primarily, this method has simple steps with a
(relatively) simple distribution which can be written explicitly. However,
the number of steps it takes to reach a given distance from the root is not
fixed and has to be estimated.

The idea behind this construction and the heuristics for the Hausdorff
dimension of the scaling limit of 2-dimensional quantum gravity may be
found in section 4.7 of \cite{ADJ}. Following them, we call the process
``peeling'', as it is similar to peeling an apple by cutting a thin strip
going around the apple in circles. The name is especially appropriate when
the peeling is made in an ordered manner, as in the Sections regarding
growth estimates, however, we use the name also in the context of
Section~\ref{sec:percolation} where the process advances in a more chaotic
manner.

Consider first the case of a free triangulation $T$ of an $(m+2)$-gon. Call
the boundary vertices $x_0,\ldots,x_{m+1}$. There is a single triangle
$t \in T$ containing the edge $(x_{m+1},x_0)$. (This edge is chosen for
simplicity; by symmetry the following discussion holds for any other
boundary edge as well). Denote the third vertex of $t$ by $y$. There are
two possibilities: either $y$ is an internal vertex of the triangulation,
or else $y=x_i$ for some $x_i$. In the former case, the rest of the
triangulation is a triangulation of an $(m+3)$-gon, hence the sum of the
weights of all such triangulations is $\alpha^{-1}Z_{m+1}$ (the factor
$\alpha^{-1}$ accounts for $y$). Thus: 
\begin{equation}\label{eq:finnew}
\mu_m(y \not \in \{x_0,\ldots,x_{m+1}\}) = \frac{Z_{m+1}}{\alpha Z_m} .
\end{equation}

\begin{figure}
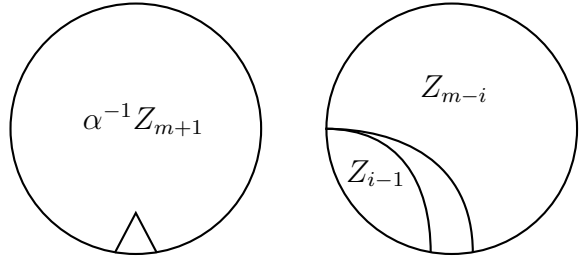

\begin{center}
\psset{unit=1.4mm}
\pspicture(-12,-12)(42,12)
  \pscircle(0,0){12}
  \psline(-2,-11.8)(0,-8)(2,-11.8)
  \put(-5,0){$\alpha^{-1} Z_{m+1}$}

  \pscircle(30,0){12}
  \psbezier(28,-11.8)(28,-5)(25,0)(18,0)
  \psbezier(32,-11.8)(32,-3)(25,0)(18,0)
  \put(27,3){$Z_{m-i}$}
  \put(20,-5){$Z_{i-1}$}
\endpspicture
\end{center}
\caption{\label{fig:free_peel}
The two possibilities when a triangle is added in a free triangulation.}
\end{figure}

Similarly, if $y=x_i$, then we have the triangle $(x_{m+1},x_0,x_i)$ and two
triangulations of an $(i+1)$-gon and an $(m-i+2)$-gon. The weight of a
triangulation, $\alpha^{-|T|}$, is multiplicative, and since any pair of
triangulations for the two components is possible, the total weight is
$Z_{i-1}Z_{m-i}$ (see Figure~\ref{fig:free_peel}). Therefore
\begin{equation}\label{eq:finold}
\mu_m(y=x_i) = \frac{Z_{i-1} Z_{m-i}}{Z_m},
\end{equation}
and conditioned on $y=x_i$ the two components are filled with independent
free triangulations.

The case $y=x_1$ is especially interesting, for then we have the triangle
$(x_{m+1},x_0,x_1)$. If this is also the triangle supported on the boundary
edge $(x_0,x_1)$, then after adding it the vertex $x_0$ is no longer on the
boundary of the remaining triangulation. Otherwise we still need to
triangulate the 2-gon $(x_0,x_1)$. This was accounted for in the above
formula since we set $\phi_{0,0}=1$ and not 0. The triangulation of a 2-gon
with no internal vertices has no internal triangles and so the edges are
glued together.

The same thing may happen in the case $y=x_m$ with the edge $(x_m,x_{m+1})$.
In the case that $m=1$ and $y=x_1=x_m$ there is a possibility that both
2-gons are glued, in which case there are no discs left to triangulate and
hence no additional triangles. 

Equations~\ref{eq:finnew} and \ref{eq:finold} allow us to sample a free
triangulation of any disc by adding triangles one at a time. Since we know
that the free triangulation is finite, the process a.s.\ terminates at some
time. Note that we can choose which boundary edge to build on at each 
iteration in any way we want to without effecting the resulting
distribution.

\medskip

A similar construction allows us to sample the UIPT. Suppose we wish to
sample the UIPT of an $(m+2)$-gon, with external boundary vertices
$x_0,\ldots,x_{m+1}$. Heuristically, the same notion of the total weight
for each possibility appears as before, with $C_m$ being the total weight
of infinite triangulations of an $(m+2)$-gon. This is so since the weight
for triangulations of size $N$ is roughly $C_m N^{-5/2}$, and the power
term will cancel with an identical term in the denominator when
$N \to \infty$.

To be precise, we consider as before the triangle $(x_{m+1},x_0,y)$
containing the edge $(x_{m+1},x_0)$. The probability that $y$ is a new
vertex is:
\begin{eqnarray}
\P(y \not \in \{x_0,\ldots,x_{m+1}\}) 
   & = & \lim_{N \to \infty} \frac{\phi_{m+1,N-1}}{\phi_{m,N}}   \nonumber\\
   & = & \frac{C_{m+1}}{\alpha C_m} .     \label{eq:infnew}
\end{eqnarray}

\begin{figure}
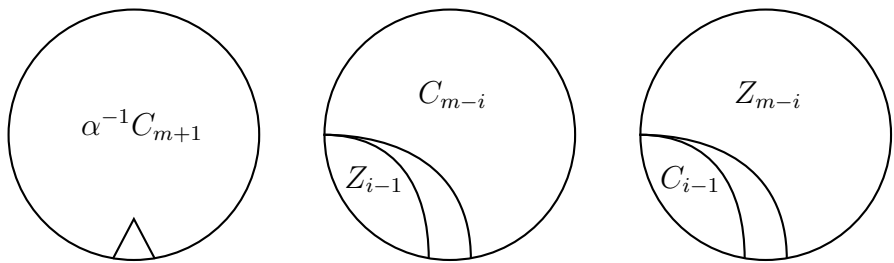

\begin{center}
\psset{unit=1.4mm}
\pspicture(-12,-12)(72,12)
  \pscircle(0,0){12}
  \psline(-2,-11.8)(0,-8)(2,-11.8)
  \put(-5,0){$\alpha^{-1} C_{m+1}$}

  \pscircle(30,0){12}
  \psbezier(28,-11.8)(28,-5)(25,0)(18,0)
  \psbezier(32,-11.8)(32,-3)(25,0)(18,0)
  \put(27,3){$C_{m-i}$}
  \put(20,-5){$Z_{i-1}$}

  \pscircle(60,0){12}
  \psbezier(58,-11.8)(58,-5)(55,0)(48,0)
  \psbezier(62,-11.8)(62,-3)(55,0)(48,0)
  \put(57,3){$Z_{m-i}$}
  \put(50,-5){$C_{i-1}$}
\endpspicture
\end{center}
\caption{\label{fig:inf_peel}
The possibilities when a triangle is added in the UIPT.}
\end{figure}

When $y=x_i$ we have two sub-triangulations: $T_1$ with boundary size $i+1$
and $T_2$ with boundary size $m-i+2$. One of them is infinite, and by
Theorem~\ref{thm:locality} the infinite one is a UIPT and the finite is a
free triangulation. We have two more possibilities, as in
Figure~\ref{fig:inf_peel}, with probabilities:
\begin{eqnarray}
\P((y=x_i) \cap T_1 \textrm{ is infinite})
   & = & \frac{Z_{m-i} C_{i-1}}{C_m}     \nonumber\\
\P((y=x_i) \cap T_2 \textrm{ is infinite})
   & = & \frac{C_{m-i} Z_{i-1}}{C_m} .  \label{eq:infold}
\end{eqnarray}

We can use (\ref{eq:infnew}), (\ref{eq:infold}) to sample a neighborhood of
the root in the UIPT, as follows: Start with the root triangle, and proceed
to add new triangles to the triangulation. When the added triangle
partitions the triangulation to a finite and infinite part, fill the finite
part with a free triangulation. In either case the remaining triangulation
is now a UIPT of some polygon.

To sample a ball of radius $r$, proceed as above as long as there are any
vertices on the boundary that are at distance less than $r$ from the root.
Since the ball is a.s.\ finite (from the tightness of its size,
Corollary~4.5 of \cite{UIPT1}), this process will a.s.\ terminate, giving
a sample of the ball of radius $r$. Termination can also be deduced
directly from the stated process by induction on $r$: Suppose that a.s.\
at some time all vertices at distance less than $r$ from the root are in
the interior, then the boundary contains only vertices at distance $r$ or
greater, so any new vertex we add is at distance at least $r+1$. Hence the
set of boundary vertices at distance $r$ can only decrease. If we extend
our triangulation at the boundary edge $(u,v)$, then with probability
bounded away from 0, $v$ will disappear from the boundary. Thus a.s.\ after
some finite number of iterations all vertices at distance $r$ will also be
in the interior. Since we have shown a process that a.s.\ terminates in
a finite time and outputs a neighborhood of the root in the UIPT, this
gives an alternative proof of the existence of the UIPT probability
measure.

\medskip

\note In the construction above there is complete freedom in the choice at
each step of which edge to build on. Typically, especially for analyzing
the growth rate of the UIPT, we will go around the triangulation in a fixed
direction. This means that we take a vertex $v$, on the boundary, and add
the triangle incident on the edge to its right as long as $v$ remains on
the outer boundary. Whenever a hole is formed we fill it with a free
triangulation. As soon as $v$ leaves the boundary we look along the old
boundary counterclockwise from $v$ to the first vertex that is in the new
boundary. This will guarantee that we find all vertices at distance $r$
from the root before proceeding to $r+1$.

For other uses we may choose differently the edge to build on. Thus for
analyzing percolation the edge we choose will depend on the colors of
vertices previously visited. The results of the next section do not depend
on the choice of the edge at each step.

\medskip

The peeling method is somewhat flexible, and lends itself to sampling from
a number of classes of planar objects. For example, in order to sample a
type I triangulation, where looped edges are allowed, one only needs to
find the appropriate values of $Z_m,C_m$ and to include the possibility
that the third vertex of the added triangle coincides with one of the other
two. If the values of $C_m,Z_m$ have the same asymptotics, then the peeling
process will proceed in much the same way and the following results may
apply to that class as well.

\begin{figure}
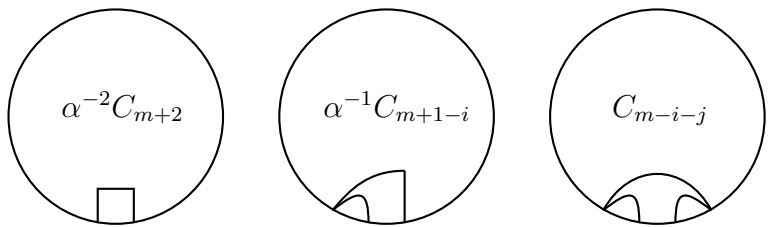

\begin{center}
\psset{unit=1.2mm}
\pspicture(-12,-12)(72,12)
  \pscircle(0,0){12}
  \psline(-2,-11.8)(-2,-8)(2,-8)(2,-11.8)
  \put(-6,0){$\alpha^{-2} C_{m+2}$}

  \pscircle(30,0){12}
  \psbezier(28,-11.8)(28,-8)(27,-8)(24,-10.392)
  \psbezier(32,-6)(28,-6)(26,-8)(24,-10.392)
  \psline(32,-11.8)(32,-6)
  \put(23,0){$\alpha^{-1} C_{m+1-i}$}

  \pscircle(60,0){12}
  \psbezier(58,-11.8)(58,-8)(57,-8)(54,-10.392)
  \psbezier(62,-11.8)(62,-8)(63,-8)(66,-10.392)
  \psbezier(66,-10.392)(63,-5)(57,-5)(54,-10.392)
  \put(55,0){$C_{m-i-j}$}
\endpspicture
\end{center}
\caption{\label{fig:quad_peel}
Some possibilities when sampling an infinite quadrangulation of the plane.}
\end{figure}

To sample a uniform infinite quadrangulation of the plane, a structure
which may be defined and proved to exist in much the same way as the UIPT,
the peeling process adds at each time a quadrangle. The number of
possibilities now grows as there may be two, one or no new vertices, as in
Figure~\ref{fig:quad_peel}. In each case also the infinite face needs to be
distinguished from the others. However, the process stays essentially the
same, and using formulas for the asymptotic number of quadrangulations of a
disc and for the partition function of the same, the process may be
analyzed similarly.

%%%%%%%%%%%%%%%%%%%%%%
\section{Markov Chain}  \label{sec:markov}
%%%%%%%%%%%%%%%%%%%%%%

To study the growth rate of the UIPT we will analyze the peeling
process for constructing the UIPT. At first we will focus only on the
evolution of the outer boundary of the generated sub-triangulations as
triangles are added.

Suppose we have a finite sub-triangulation with outer boundary of size
$m+2$. We add to it the triangle on the outer face containing a given
boundary edge $e$. Two of the vertices of this triangle are given by the
ends of $e$. The probability that the third vertex is not in the
sub-triangulation is given by (\ref{eq:infnew}), and on this event the size
of the boundary increases by 1. The probability that the third vertex is
some vertex on the outer boundary is given by (\ref{eq:infold}). In this
case the triangulation of the outer face is partitioned to a finite part
with boundary of size $i+1$ for some $i \ge 1$, and an infinite part with
boundary of size $m-i+2$. Thus the size of the outer boundary decreases by
$i$. The probability of that is $2\frac{C_{m-i}Z_{i-1}}{C_m}$. (The factor
of 2 comes from the fact that there are 2 vertices at distance $i$ from
$e$).

Let $M_n$ be the size of the outer boundary after $n$ such triangles have
been added. We see that the evolution of $M_n$ can be described as a Markov
chain satisfying $M_{n+1} = M_n + X_n$ where the distribution of $X_n$ is
given by:
\begin{eqnarray}
\P(X_n=1|M_n=m)
   &=& \frac{C_{m+1}}{\alpha C_m} = \frac{2m+3}{3m+3} ,   \nonumber  \\
\P(X_n=-k|M_n=m)
   &=& \frac{2C_{m-k}Z_{k-1}}{C_m}      \label{eq:step_dist} \\
   &=& \frac{2(2k-2)!}{(k-1)!(k+1)!} \cdot
       \frac{m!^2(2m-2k+1)!}{(m-k)!^2(2m+1)!}  \nonumber  
\end{eqnarray}
for $k \in \N$.

Denote these probabilities by $p_{-1,m}$ and $p_{k,m}$ respectively.
When $m$ is large these probabilities converge to a limit distribution.
Denote by $p_{-1}$ (resp. $p_k$) the limit of the probability of having
$X=1$ (resp. $X=-k$), then:
\begin{eqnarray}
p_{-1} &=& \lim_{m \to \infty} p_{-1,m} 
   =  \frac{2}{3} ,                   \nonumber  \\
p_k &=& \lim_{m \to \infty} p_{k,m} 
   =  \frac{2(2k-2)!}{(k-1)!(k+1)!4^k} .  \label{eq:limit_dist}
\end{eqnarray}

\note It is worthwhile observing the relation:
\[
p_{k,m} = p_k \frac{(m)_k}{(m+1/2)_k} ,
\]
where
\[
(x)_k = x(x-1)\ldots(x-k+1) = \frac{\Gamma(x+1)}{\Gamma(X-k+1)}
\]
is the descending factorial notation. This implies that for $k>0$, the
probabilities $p_{k,m}$ are increasing in $m$ and converge to $p_k$. From
this it follows that if $m_1<m_2$ then the step distribution from $m_1$
(i.e.\ the law of $X_n$ conditioned on $M_n=m_1$) stochastically dominates
the step distribution from $m_2$. 

Another consequence is that for any $a<1$ there is a $c=c(a)>0$, such that
$k<am$ implies $p_{k,m}>cp_k$. I.e., the probability of steps of size up to
a constant fraction of $m$ is within a constant factor of the limit
probability. Of course for steps larger than $m$ the probability is 0.
Hence the step distribution from $m$ is similar to the limit distribution
conditioned to be at most $m$.

\medskip

In light of the convergence of the step distributions, the states of the
Markov chain, $M_n$, can be viewed as a random walk on the integers, with
step distribution depending slightly on the location. The Markov chain will 
be shown to be transient, and so the step distributions indeed converge. A
central difficulty arises from the fact that the limit distribution has 
infinite variation (indeed $\textrm{Var}(X|M=m) \sim c\sqrt m$).

Let $X$ be a sample of the limit distribution. Because $p_k \sim ck^{-5/2}$,
$X$ has a finite $\alpha$'th moment iff $\alpha<3/2$. Since $\E X=0$ and
$X$ is bounded by 1, the theory of stable random variables (see
\cite{SaTa,Zolo}) tells us that if $X_0,X_1,\ldots$ were i.i.d.\ copies of
$X$, then $n^{-2/3} \sum_{i<n} X_i$ would have converged in distribution to
a totally asymmetric stable random variable of type 3/2. Of course, in our
case, the steps $X_i$ are neither independent nor are they equally
distributed. Instead their distribution depends on the sum of their
predecessors. Still, this gives an indication that $M_n$ should be studied
at the scale of $n^{2/3}$. In fact, the sequence $M_n$ with proper scaling
appears to converge to a stable process conditioned to remain positive.

\medskip

Since $M_n$ is finite, the expectation of $X_n$ is not 0, but rather the
random walk has some drift. A straightforward calculation using generating
functions leads to the following result:

\begin{lemma}\label{lem:step_exp}
\[
\E(X_n|M_n=m) = \frac{4^m m!^2}{(2m+1)!}
        \sim \sqrt{\frac{\pi}{2m}} .
\]
\end{lemma}

\begin{proof}[Sketch of proof]
In view of the distribution of $X_n$, calculating the expectation involves
substituting $x=1$ and finding the coefficient of $y^m$ in
\[
\sum_{k,m} y^m C_{m-k} Z_{k-1} k x^k .
\]
This may be re-written as
\[
\sum_k k Z_{k-1} (xy)^k \sum_{m\ge k} C_{m-k} y^{m-k} .
\]
However, using the binomial formula we have:
\[
\sum C_l y^l = \frac{\sqrt{3}}{4\sqrt{\pi}} (1-9y)^{-3/2} ,
\]
and
\[
\sum k Z_{k-1} (xy)^k = \frac{2-(2+9x)\sqrt{1-9x}}{54x} .
\]
Giving a closed form for the above double sum.
\end{proof}

This leads us to believe that the rate of growth of $M_n$ is roughly equal
to $M_n^{-1/2}$. This corresponds (again) to a growth rate of $n^{2/3}$. In
fact, a lower bound on $\E M_n$ follows from the convexity of
$x \to x^{-1/2}$: 
\[
\E X_n \ge c \E (M_n^{-1/2}) \ge c (\E M_n)^{-1/2} ,
\]
\[
\E M_{n+1} = \E(M_n+X_n) \ge \E M_n + c (\E M_n)^{-1/2} ,
\]
and therefore $\E M_n \ge c' n^{2/3}$.

This rate of growth is indeed correct and we prove the following variation
on the LIL for $M_n$:

\begin{thm}\label{thm:M_small}
A.s.
\[
\limsup \frac{M_n}{n^{2/3} \log n} < \infty .
\]
\end{thm}

The proof is similar to the proof of the law of iterated logarithms, with
some modifications to accommodate the positive expectation and unbounded
variation of the steps. While a better upper bound with an iterated
logarithm may hold, we will not attempt to prove such an upper bound here,
since it is the power term we are primarily interested in. Additionally, in
evaluating the growth rate of the UIPT, other steps add logarithmic
factors, so the end result will not improve significantly by having a
tighter bound here.

Let $\F_n$ denote the $\sigma$-field generated by the random variables
$M_0,\ldots,M_n$. Thus $\F_0$ is the trivial $\sigma$-field. $\E_{\F_n}$
will denote expectation w.r.t.\ $\F_n$, i.e., expectation conditioned on
the past. Lemma~\ref{lem:step_exp} may therefore be stated as
$\E_{\F_n} X_n \sim cM_n^{-1/2}$.

\begin{proof}
We use the fact that for any $\lambda_0$ there is an absolute constant
$c_1$, such that for $t < \lambda_0$
\[
e^t < 1 + t + c_1 |t|^{3/2} .
\]

Denote by $\barX_n$ the normalized steps: $\barX_n = X_n - \E_{\F_n}X_n$.
Using the bounds $0<\E_{\F_n} X_n\le 1$, the $3/2$'th moment of $\barX_n$
is estimated.
\begin{eqnarray*}
\E_{\F_n} |\barX_n|^{3/2} 
  &=& p_{-1,M_n}(1-\E_{\F_n}X_n)^{3/2} +
      \sum_{k=1}^{M_n} p_{k,M_n} (k+\E_{\F_n}X_n)^{3/2}     \\
  &<& 1 + \sum_{k=1}^{M_n} p_k (k+1)^{3/2}      \\
  &<& c_2 \sum^{M_n} k^{-1}                     \\
  &<& c_3 \log M_n                              \\
  &<& c_3 \log n .
\end{eqnarray*}
A lower bound with a different constant also holds hence this is best
possible.

For any small $\lambda>0$ we have
\begin{eqnarray*}
\E_{\F_n} e^{\lambda \barX_n}
  & \le &  \E \left( 1+\lambda \barX_n + 
               c_1 (\lambda |\barX_n|)^{3/2} \right)    \\
  & \le &  1 + c_4 \lambda^{3/2} \log n             \\
  &  <   &  \pexp{c_4 \lambda^{3/2}\log n} .
\end{eqnarray*}

While $\barX_n$ are not independent, the last bound is deterministic
(independent of the history). Therefore we may multiply for $m\le i < n$
to get:
\[
\E \pexp {\lambda \sum_{i=m}^{n-1} \barX_i} <  
   \pexp {c_4 \lambda^{3/2} (n-m) \log n} ,
\]
and so:
\[
\P\left(\sum_{i=m}^{n-1} \barX_i > t\right )
  < \pexp{-t\lambda + c_4 \lambda^{3/2} (n-m) \log n} .
\]
For small enough $\lambda$, this is a bound on the probability that the sum
of the steps over an interval $[m,n)$ is much bigger than the sum of their
expectations. $\lambda$ will be chosen to minimize this bound w.r.t.\
$t,m,n$. The optimal $\lambda$ is such that
\[
\sqrt \lambda = 2c_4 \frac{t}{(n-m)\log n} ,
\]
and then the bound becomes
\[
\P\left(\sum_{i=m}^{n-1} \barX_i > t\right )
  <  \pexp{-c_5 \frac{t^3}{(n-m)^2 \log^2 n}} .
\]
If $t$ is such that $t^3 = (3/c_5) (n-m)^2 \log^3 n$, then
$\lambda  = O((n-m)^{-2/3})$ tends to 0 as the size of the interval is
large. Substituting the given values for $t$ and $\lambda$ we get for
sufficiently large intervals:
\[
\P \left ( \sum_{i=m}^{n-1} \barX_i > c_6 (n-m)^{2/3} \log n \right) 
   < n^{-3} .
\]
For small intervals this can be made to hold as well by increasing $c_6$.

Since the sum over all intervals $[m,n)$ with $m<n$ of $n^{-3}$ is finite,
by Borel-Cantelli, all but finitely many intervals satisfy:
\begin{equation}\label{eq:goodwin}
\sum_{i=m}^{n-1} \barX_i  < c_6 (n-m)^{2/3} \log n .
\end{equation}

We use now an upper bound on $\E_{\F_n} X_n$: There is an $a>0$ such that
$\E_{\F_n} X_n < aM_n^{-1/2}$. Take $c_7=c_6+1+a$. For each $n$ where
$M_n>c_7 n^{2/3} \log n$ we will find an interval $[m,n)$ where this bound
is violated, proving the theorem.

Suppose $n$ is such that $M_n>c_7n^{2/3}\log n$. Take
\[
m = 1+\max \{ m<n | M_m \le n^{2/3}\}.
\]
By the choice of $m$, for each $i \in [m,n)$ we have $M_i>i^{2/3}$ and
therefore $\E_{\F_i} X_i < a i^{-1/3}$. We then have:
\[
\sum_{m}^{n-1} \E_{\F_i} X_i < (n-m) an^{-1/3} < an^{2/3} .
\]
On the other hand:
\begin{eqnarray*}
\sum_m^{n-1} \barX_i
  &=& M_n-M_m - \sum_m^{n-1} \E_{\F_i} X_i          \\
  &>& c_7 n^{2/3} \log n - n^{2/3} - an^{2/3}       \\
  &>& c_6 n^{2/3} \log n ,
\end{eqnarray*}
contradicting (\ref{eq:goodwin}), thus we found an interval $[m,n)$ as
claimed.
\end{proof}

Next, we consider the hitting probabilities of the Markov chain. Those will
later be used in establishing a lower bound. The proof is another
straightforward use of generating functions, similar to the proof of
Lemma~\ref{lem:step_exp}.

\begin{claim}\label{clm:hitprob}
Starting the Markov chain at $n$, the probability of it ever hitting $m$ is
given by:
\[
\P(n \to m) = 1-\frac{(n)_{m+1}}{(n+1/2)_{m+1}} .
\]
\end{claim}

\begin{proof}[Sketch of proof]
Fix $m$ and assume the probability in question is $a_n$, with $a_m=1$. For
any $n \neq m$ we must have
\begin{eqnarray*}
a_{n,m} &=& \sum_{k=-1}{n} p_{n,k} a_{n-k,m}    \\
        &=& \frac{C_{n+1}}{\alpha C_n} a_{n+1,m} +
            \sum_{k=1}{n} \frac{2C_{n-k}Z_{k-1}}{C_n} a_{n-k,m} .
\end{eqnarray*}
If $g(t)=(t\alpha)^{-1}-1+2\sum Z_{k-1}t^k$, and $f(t)=\sum C_na_nt^n$,
then this just says that
\[
f(t)g(t) = \frac{C_0}{\alpha t} + \beta t^m ,
\]
where $\beta$ is determined by boundary conditions. Since
$g(t)=2(1-9t)^{3/2}/27t$ is known, we can find $f(t)$, express as a power
series and divide by $C_n$ to get $a_n$. 
\end{proof}

\note In particular, the probability of returning to 0 from $n$ is
$\frac{1}{2n+1}$ and thus the Markov chain is transient. This formula also
holds for $m\ge n$, since the numerator vanishes. For $n \gg m$ the hitting
probability can is $\frac{m+1}{2n} + O(n^{-2})$.

As a simple consequence, the probability of never returning to $n$ when
starting from $n$ is found. The only way to avoid returning to $n$ is by
having $X=1$ and afterward not hitting $n$ from $n+1$. Since every number
is visited at least once, we get:

\begin{claim}\label{clm:localtime}
The number of visits to $n$ when the Markov chain is started from 0, is a
geometric random variable with mean 
\[
\frac{3n+3}{2n+3} \frac{(n+3/2)_{n+1}}{(n+1)!}  \sim c\sqrt n .
\]
\end{claim}

\begin{thm}\label{thm:M_large}
For any $\ep>0$, a.s.
\[
\lim_{n \to \infty} \frac{M_n \log^{2+\ep} n}{n^{2/3}} = \infty .
\]
\end{thm}

\begin{proof}
Consider the time $T_n$ when the Markov chain first hits $2^n$. Fix some
$\ep>0$ and set $a_n=2^n n^{-(1+\ep)}$. We consider the behavior of the
Markov chain from time $T_n$ until it leaves the interval
$I_n = [a_n, 2^{n+1}]$. For this we consider the probability of the Markov
chain taking a step from $I_n$ into the interval $J_n=[a_n/2,a_n]$. The
probability of taking a step from $M$ to $M'<M$ is monotone increasing in
$M'$ and monotone decreasing in $M$, hence among all steps from $M\in I_n$
to $M'\in J_n$ the step from $2^{n+1}$ to $\lceil a_n/2 \rceil$ has the
smallest probability. Using Stirling's formula, the probability of such a
step is 
\begin{eqnarray*}
\P\left(M' = a_n/2 | M=2^{n+1} \right)
  & = & \P(X = a_n/2 - 2^{n+1} | M=2^{n+1})       \\
  & > & c_1 \frac{(2^n)^{-5/2}}{n^{(1+\ep)/2}} .
\end{eqnarray*}
It follows that the probability of taking a step into $J_n$ from any
$M \in I_n$ is at least
\[
c_1 \frac{(2^n)^{-5/2}}{n^{(1+\ep)/2}} |J_n| =
c_2 \frac{(2^n)^{-3/2}}{n^{3(1+\ep)/2}} .
\]
Therefore for any $k$, 
\begin{eqnarray*}
\P(\textrm{visiting $J_n$ after $T_n$})
  &\ge& \P(T_{n+1}-T_n > k) \cdot  \\
  &   & \qquad \P(\textrm{visiting $J_n$ after $T_n$}|T_{n+1}-T_n>k) \\
%        \P(T_{n+1}-T_n > k)                                     \\
  &\ge& \left(1-\left(1-c_2 \frac{(2^n)^{-3/2}}{n^{3(1+\ep)/2}}
        \right)^k\right)  \P(T_{n+1}-T_n > k) .
\end{eqnarray*}
Setting $k=n^{3(1+\ep)/2} 2^{3n/2}$ so that the first factor is roughly
constant, gives
\[
\P(T_{n+1}-T_n > n^{3(1+\ep)/2} 2^{3n/2}) \le
c_3 \P(\textrm{visiting $J_n$ after $T_n$}) .
\]

However, by Claim~\ref{clm:hitprob} the probability that after time $T_n$
the Markov chain visits $\lfloor a_n \rfloor$ is approximately
$\frac{1}{2n^{1+\ep}}$. Since the Markov chain can increase by at most 1 at
each step, if $J_n$ is visited after time $T_n$, then so is
$\lfloor a_n \rfloor$. Therefore
\[
\P(T_{n+1}-T_n > n^{3(1+\ep)/2} 2^{3n/2}) \le \frac{c_4}{n^{-(1+\ep)}} .
\]
Thus a.s.\ $T_{n+1}-T_n < n^{3(1+\ep)/2} 2^{3n/2}$ for all but finitely many
$n$. Summing, we get that for some $c_5$, a.s.\ for all but finitely many
$n$, $T_n < c_5 n^{3(1+\ep)/2} (2^n)^{3/2}$.

This gives a sequence of points lying below the graph of $M_n$. By
interpolating this translates to the lower bound on $M_n$ for any $n$.
Formally, let $m$ be such that $T_m \le n < T_{m+1}$. Since the Markov
chain reached $2^m$ before time $n$, necessarily $n\ge 2^m$, and by the 
bound on $T_{m+1}$, for all but finitely many $n$: 
\begin{eqnarray*}
n < T_{m+1}  &<&  c_5 (m+1)^{3(1+\ep)/2} (2^{m+1})^{3/2} ,    \\
    n^{2/3}  &<&  c_6 m^{1+\ep} 2^m .
\end{eqnarray*}
Since $J_m$ is not visited after $T_m$ (for large $m$),
\begin{eqnarray*}
M_n &>& 2^m m^{-(1+\ep)}          \\
    &>& c_6^{-1} n^{2/3} m^{-2(1+\ep)}      \\
    &>& c_7 n^{2/3} \log^{-2(1+\ep)} n .
\end{eqnarray*}
\end{proof}

%%%%%%%%%%%%%%%%%%%%%%%%%
\section{Boundary Growth}  \label{sec:boundary}
%%%%%%%%%%%%%%%%%%%%%%%%%

Based on Theorems~\ref{thm:M_small} and \ref{thm:M_large} we now turn to
study the growth rate of the UIPT. The UIPT is sampled starting with only
a root triangle, and adding triangles in the outer face of the
sub-triangulation we have. The behavior of $M_n$ --- the size of the outer
boundary after $n$ triangles have been added --- is known. We wish now to 
relate this to the distance between the root and the boundary.

Throughout this and the next two sections triangles are added in an ordered
manner, going around the boundary counterclockwise, and adding all the
triangles incident on a vertex before moving to the next. In this manner we
are sure to encounter all vertices at distance $r$ from the root before
moving on to $r+1$.

\medskip

Initially, we investigate the number of steps it takes the above process to
sample the ball of radius $r$ around the root (as well as its hull). This
translates to estimates on the size of the ball's boundary. Let $T_r$ be
the time (number of triangles added) when we have found the outer boundary
of the ball of radius $r$ around the root. We prove:

\begin{thm}\label{thm:layertime}
A.s.
\[
\limsup_{r \to \infty} \frac{T_r}{r^3 \log^3 r} < \infty ,
\]
and for any $\ep>0$
\[
\lim_{r \to \infty}  \frac{T_r \log^{6+\ep} r}{r^3} = \infty .
\]
\end{thm}

Let $C_n$ denote the set of vertices on the outer boundary at time $n$ (thus
$|C_n| = M_n+2$). Clearly after some number of steps we have gone round a
full circle, and added every triangle incident on a vertex in $C_n$. Denote
the time it takes to accomplish this by $D_n$, i.e.,
\[
D_n = \inf\{d|C_{n+d} \cap C_n = \emptyset\} .
\]

Clearly $T_{r+1} = T_r + D_{T_r}$. It is plausible that the process
proceeds in an approximately fixed rate along the boundary, so that $D_n$
is roughly linear in $M_n$. We therefore expect $T_{r+1}-T_r=D_{T_r}$ to
be on the order of $M_{T_r} \approx T_r^{2/3}$, indicating that $T_r$ grows
like $r^3$. Together with Theorems~\ref{thm:M_small} and \ref{thm:M_large}
this implies that $|C_{T_r}|$ grows like $r^2$, i.e., the size of the outer
boundary of the ball of radius $r$ is quadratic in $r$. Note that $T_r$ is
not the volume of the ball but only the number of steps made. At each step
when $X_t<0$, a free triangulation is glued to the sub-triangulation. The
bulk of the volume lies in those.

We prove first the following estimate:

\begin{lemma} \label{lem:cycletime}
Starting with $M$ boundary vertices, let $D$ denote the number of steps of 
the peeling process until every vertex originally on the boundary is in the
interior. For some $a,c>0$, $\P(D>a M) < e^{-cM}$ .
\end{lemma}

\note It is impossible to get an exponential bound for $\P(D \le bM)$ for
any $b$, due to the thick tail on the negative side of $X_n$. In fact, for
every $b>0$ there will a.s.\ be infinitely many $n$'s for which $D_n<bM_n$.
This is so, since the expected number of visits to $m$ is approximately
$\sqrt m$ (by Claim~\ref{clm:localtime}) and the probability of making a
single step of size at least $(1-\ep)m$ is about $m^{-3/2}$, thus such
large jumps are made infinitely often.

\begin{proof}[Proof of Lemma~\ref{lem:cycletime}]
Let $Z_t$ be the number of vertices in $C_0$ that are removed from the
boundary at time $t$. $D$ is just the smallest number such that
$\sum^D Z_t = M$, since once a vertex is in the interior of the
triangulation it can not return to the boundary at a later stage. When
building a triangle on an edge, with probability $1-p_{-1,M_t}$ the new
triangle does not introduce a new vertex, and some vertices leave the outer
boundary. Since an endpoint of the edge is in $C_0$, with half that
probability the vertices leaving the boundary are from $C_0$. Hence at each
step, $\P(Z_t>0)$ is bounded away from 0.

Since vertices are removed from the boundary at a positive rate bounded
from 0, for some $a$ the probability of taking more than $a M$ steps for
all vertices in $C_0$ to be removed decays exponentially in $M$.
\end{proof}

\begin{proof}[Proof of Theorem~\ref{thm:layertime}]
By Theorem~\ref{thm:M_large}, $M_{T_r}$ grows like a power of $T_r$. By
Lemma~\ref{lem:cycletime}, the probability that $D_{T_r} > aM_{T_r}$ decays
exponentially in $M_{T_r}$. It follows that a.s.\ for all but finitely many
$r$ we have
\begin{eqnarray*}
D_{T_r} & < & a M_{T_r} < c_1 T_r^{2/3}\log T_r ,       \\
T_{r+1} & < & T_r + c_1 T_r^{2/3} \log T_r ,
\end{eqnarray*}
and so
\[
T_r < c_2 r^3 \log^3 r .
\]

\medskip

For the lower bound, note that the vertices are added to the boundary only
one at a time (when $X_n=1$). Since all the vertices in $C_{T_{r+1}}$ have
been added to the boundary after time $T_r$, it follows that
\[
T_{r+1} - T_r > M_{T_{r+1}}.
\]
Together with Theorem~\ref{thm:M_small} this implies that for any $\ep>0$,
a.s.\ for all but but finitely many $r$,
\[
T_{r+1} - T_r > T_{r+1}^{2/3} \log^{-(2+\ep)} T_{r+1}
              > T_r^{2/3} \log^{-(2+\ep)} T_r.
\]
And so for any $\ep>0$, for large $r$,
\[
T_r > r^3 \log^{-(6+\ep)} r.
\]
\end{proof}

Combining this Theorem with Theorems~\ref{thm:M_small} and \ref{thm:M_large}
gives:

\begin{coro} \label{:cor:boundary_size}
The size of the boundary of the $B_r$'s hull, given by $M_{T_r}$, a.s.\
satisfies
\[
\limsup_{r \to \infty} \frac{M_{T_r}}{r^2 \log^3 r} < \infty ,
\]
and for any $\ep>0$
\[
\lim_{r \to \infty} \frac{M_{T_r} \log^{6+\ep} r}{r^2} = \infty .
\]
\end{coro}

%%%%%%%%%%%%%%%%%%%%%%%%%%%%
\section{Hull Volume Growth}  \label{sec:hull_vol}
%%%%%%%%%%%%%%%%%%%%%%%%%%%%

Knowing how long it takes to sample $B_r$ in the UIPT, we turn our
attention to the distribution of its hull's volume. So far, new vertices
were added at times when $X_t=1$. Since $\P(X_t=1)$ remains away from 0,
the number of vertices added that way is linear in $T_r$. This shows a
volume growth of at least $r^3$ (up to logarithmic factors).

There are many additional vertices that are added whenever $X_t<0$. At
those times a portion of the boundary is closed off and filled with a free
triangulation of an $(|X_t|+1)$-gon. Using Proposition~\ref{prop:Z} and the
derivative of $Z_m(t)$ we find:

\begin{prop}\label{prop:free_size} 
The expected number of internal vertices in a free triangulation of an
$(m+2)$-gon is:
\[
(Z_m)^{-1} \sum n \phi_{n,m} \alpha^{-n} = \frac{(m+1)(2m+1)}{3} .
\]
\end{prop}

The heuristic explanation for the volume growth is as follows: Let $Y_t$ be
the number of vertices added at time $t$. The idea behind the proof of the
volume growth is as follows: By the above, $\E(Y_t|X_t) \sim c X_t^2$.
Since $\P(X_t=-k) \approx ck^{-5/2}$ for $k<M_t$, and the probability is 0
for $k>M_t$, the expected number of vertices added at time $t$ is (up to
multiplicative constants): 
\[
\E(Y_t|M_t) \approx \sum \P(X_t=-k) \E(Y_t|K_t=-k)
       \approx \sum^{M_t} ck^{-1/2}
       \approx cM_t^{1/2} .
\]

Now, $M_t \approx t^{2/3}$ and therefore the expected number of vertices
added at time $t$ is roughly $t^{1/3}$. Summing that up to time $T_r$, the
number of vertices added is roughly $T_r^{4/3} \approx r^4$. This implies a
growth of $r^4$ as was suggested in the physics literature (\cite{AmWa,ADJ}
and others).

Note that in the above discussion it was convenient to consider instead of
$B_r$, its hull. This enables us to add a free triangulation which may
include vertices further away from the root. To estimate the ball volume
growth the distance of vertices in such free triangulations need to be
considered as well.

For the steps $X_t$ of the Markov chain with distribution given by
(\ref{eq:step_dist}), define
\begin{eqnarray*}
V_T(\gamma) &=& \sum^T |X_t|^\gamma .   \\
\end{eqnarray*}

\begin{lemma} \label{lem:X_small}
Let $V_T(\gamma)$ be as above. For $\gamma>3/2$, and for any $\ep>0$, a.s.
\[
\limsup \frac{V_T(\gamma)}{T^{2\gamma/3} \log^{2\gamma/3+\ep} T} < \infty.
\]
\end{lemma}

\begin{proof}
Since the distribution of $|X_t|$ is stochastically increasing w.r.t.\
$M_t$, and since the limit distribution (\ref{eq:limit_dist}) satisfies
$\P(|X|>\lambda) = O(\lambda^{-3/2})$,
\[
\P \big(|X_t|^\gamma > \lambda|M_t\big) = O(\lambda^{-3/(2\gamma)}) 
\]
uniformly for all $t$. It follows that for some constants $a,b>0$
(possibly depending on $\gamma$) we have stochastic domination by a stable
random variable $|X_t|^\gamma \slt aS_{3/(2\gamma)}+b$. Since this
domination holds even when conditioning on $M_t$, the sum is dominated by
i.i.d.\ copies giving
\[
V_T(\gamma) = \sum^T |X_t|^\gamma \slt aT^{2\gamma/3}S_{3/(2\gamma)} + bT.
\]
Therefore as $T \to \infty$
\begin{eqnarray*}
\P\left(V_T(\gamma) > T^{2\gamma/3} \log^{2\gamma/3+\ep} T \right) 
  &<& \P\left(aS_{3/(2\gamma)} > \log^{2\gamma/3+\ep} T - o(1)\right)       \\
  &=& O(\log^{-(1+3\ep/(2\gamma))} T) .
\end{eqnarray*}

Restricting our attention to times $T_n=2^n$, we see that a.s.\ for only
finitely many such times
$V_T(\gamma) > T^{(2\gamma)/3} \log^{2\gamma/3+\ep} T$.
By the monotonicity in $T$ of $V_T(\gamma)$, the desired result holds for
all $T$.
\end{proof}

While the $\log T$ term may be replaceable by an iterated logarithm, the
exponent of $T$ is the correct one, as is evident from the following lemma.

\begin{lemma} \label{lem:X_large}
Let $X_t$ be the steps of the Markov chain with distribution described by
(\ref{eq:step_dist}). Then for any $\ep>0$, a.s.\ 
for all but finitely many $n$
\[
\big|\big\{
t \in [2^n,2^{n+1}) \textrm{ s.t.\ } |X_t|>(2^n)^{2/3} n^{-(2+\ep)}
\big\}\big| > n^2.
\]
\end{lemma}

\begin{proof}
First, estimate the probability that $|X_t|$ is very large. As long as
$\gamma < M_t/3$
\begin{eqnarray*}
\P(|X_t| \ge \gamma | M_t)
  &=& \sum_{k=\gamma}^{M_t} \P(X_t=-k|M_t)  \\
  &>& \sum_{\gamma}^{M_t/2} c_1 p_k         \\
  &>& \sum_{\gamma}^{M_t/2} c_1 k^{-5/2}    \\
  &>& c_2 \gamma^{-3/2} .
\end{eqnarray*}
Thus the probability that $|X_t|<\gamma$ for all $t \in [2^n,2^{n+1})$
satisfies
\begin{eqnarray*}
\P\left( \max_{2^n \le t <2^{n+1}} \{ |X_t| \} < \gamma \right)
  &<& \left(1-c_2\gamma^{-3/2}\right)^{2^n}                     \\
  &<& \pexp{-c_2 2^n \gamma^{-3/2}} .
\end{eqnarray*}

Set $\gamma = \gamma(n) = n^{-(2+\ep)} (2^n)^{2/3}$.
Theorem~\ref{thm:M_large} implies that a.s.\ for all but finitely many $n$
and for all $t\ge 2^n$, $M_t>3\gamma$, and hence the last bound is valid.
Thus 
\[
\P\left( \max_{2^n \le t <2^{n+1}}\{|X_t|\}<n^{-(2+\ep)}(2^n)^{2/3} \right)
  < \pexp{-c_2 n^{3+3\ep/2}} .
\]
Since this is summable, a.s.\ for all but finitely many $n$ there is a
large $|X_t|$. Since $\ep$ is arbitrary this concludes the proof.

Similarly, the probability that there are at most $k$ times with
$|X_t|<\gamma$ is bigger by at most a binomial coefficient and so is at
most 
\[
\binom{2^n}{k}
\left(1-c_2\gamma^{-3/2}\right)^{2^n-k}
  < 2^{nk} \pexp{-c_3 n^{3+3\ep/2}} .
\]
For $k=n^2$, this bound too is summable, hence a.s.\ for all but finitely
many $n$ there are at least $n^2$ large $|X_t|$'s in the interval
$[2^n,2^{n+1})$.
\end{proof}

Next, we consider the volume of the hulls. At any step in the peeling
construction where $X_t<0$, a free triangulation of a $(|X_t|+1)$-gon is
added to the triangulation. Denote by $Y_t$ the number of vertices added to
the triangulation at step $t$. If $X_t=1$, then $Y_t=1$, otherwise it is
the size of the free triangulation added.

\begin{lemma} \label{lem:Y_small}
For any $\ep>0$, a.s.
\[
\limsup_{T \to \infty} \frac{\sum^T Y_t}{T^{4/3} \log^{2+\ep} T} < \infty .
\]
\end{lemma}

The basic tool used in the proof is the following lemma:

\begin{lemma} \label{lem:prob_dom}
There are constants $c_1,c_2$ such that conditioned on $X_t=-k$, for any
$\gamma$ there is a coupling of $k^{-2} Y_t$ and  $c_1 S_{3/2} + \gamma$
such that the latter is larger with probability at least
$1-e^{-c_2\gamma}$. This is written as
\[
k^{-2} Y_t \slt c_1 S_{3/2} + \gamma
                            \qquad \textrm{with prob. $1-e^{-c_2\gamma^3}$} .
\]
\end{lemma}

\begin{proof}[Proof of Lemma~\ref{lem:Y_small}]
Lemma~\ref{lem:prob_dom} may be rephrased as: For each $t$, conditioned on
$X_t$:
\[
Y_t \slt c_1 X_t^2 S_{3/2} + \gamma X_t^2
                        \qquad \textrm{with prob. $1-e^-{c_2\gamma^3}$} .
\]
When conditioning on the whole sequence $\{X_t\}$, the $Y_t$'s become
independent. Summing up to time $T$, and using the formula for the sum of
independent stable random variables we get:
\[
\sum^T Y_t \slt c_1 V_T(3)^{2/3} S_{3/2} + \gamma V_T(2)
                        \qquad \textrm{with prob. $1-Te^{-c_2\gamma^3}$} ,
\]
and therefore 
\[
\P\left(\sum^T Y_t > \gamma V_T(2) + c_1 V_T(3)^{2/3} \log^{2/3+\ep} T
  \right) <  \P(S_{3/2} > \log^{2/3+\ep} T) + Te^{-c_2\gamma^3} .
\]
Set $\gamma = c_3 \log^{1/3} T$ with $c_3=(2/c_2)^{1/3}$. The second term
in the RHS is $T^{-1}$, and is dominated by the first. Thus
\[
\P\left(\sum^T Y_t > c_3 V_T(2) \log^{1/3} T +
                     c_1 V_T(3)^{2/3}\log^{2/3+\ep} T \right)
  < c_4 \log^{-(1+3\ep/2)} T .
\]

Considering only the times $T=2^n$, these failure probabilities are
convergent and so a.s.\ failure occurs only a finite number of those times.
Since $\sum^T Y_t$ is monotone as are $V_T(2),V_T(3)$, it follows that
a.s.\ for all but finitely many $T$
\[
\sum^T Y_t < c_3 V_{2T} \log^{1/3} T + c_1 V_{2T}(3)^{2/3} \log^{2/3+\ep} T) .
\]

Using Lemma~\ref{lem:X_small} to estimate $V_T(2),V_T(3)$ we get:
\[
\lim \frac{c_3 V_{2T} \log^{1/3} T}{T^{4/3} \log^{2+\ep} T} < \infty,
\]
as well as
\[
\limsup \frac{c_1 V_{2T}(3)^{2/3} \log^{2/3+\ep} T}{T^{4/3} \log^{2+\ep} T}
  < \infty.
\]
\end{proof}

\begin{proof}[Proof of Lemma~\ref{lem:prob_dom}]
By Theorem~\ref{thm:locality}, conditioned on $X_t=-k$, $Y_t$ is
distributed as the size of a free triangulation of a $(k+1)$-gon. From the
formula for $\phi_{n,k-1}$ given by Proposition~\ref{prop:count}, it can be
seen that
\[
\P(k^{-2}Y_t > u | X_t=-k) < c_3 u^{-3/2}
\]
for some universal $c_3$. For a real random variable $R$ use the notation 
$\Phi_R(t)=\P(R<t)$. The asymptotics of $\Phi_Y(t), \Phi_{S_{3/2}}(t)$
imply that for some universal $t_0,c_1$, for $t>t_0$
\[
\Phi_{c_1S_{3/2}}(t) < \Phi_Y(t).
\]

For $t>t_0$, monotonicity implies
\[
\Phi_{c_1S_{3/2}}(t-\gamma) < \Phi_Y(t),
\]
while for $t\le t_0$
\[
\Phi_{c_1S_{3/2}}(t-\gamma) < \Phi_{c_1S_{3/2}}(t_0-\gamma) <
  e^{-c_2\gamma^3}.
\]
Thus for any $t$
\[
\Phi_{c_1S_{3/2}}(t-\gamma) < \Phi_Y(t) + e^{-c_2\gamma^3},
\]
which is the claimed domination.
\end{proof}

The estimates resulting from comparing $Y_t$ to a stable r.v.\ $S_{3/2}$
are tight. For one thing, $Y_t$ also stochastically dominates a suitable
normalized stable random variable. A more direct proof that the exponent
4/3 of Lemma~\ref{lem:Y_small} is correct comes from the following lemma.

\begin{lemma} \label{lem:Y_large}
Let $Y_t$ be the number of vertices added to the triangulation at time $t$.
For any $\ep>0$, a.s.
\[
\lim_{T \to \infty} \frac{n^{8/3+\ep}}{(2^n)^{4/3}} 
\max_{2^n\le t<2^{n+1}} \{Y_t\} = \infty .
\]
\end{lemma}

\begin{proof}
Consider only the number of vertices added at times when $|X_t|$ is large.
Fix $\ep>0$. By Lemma~\ref{lem:X_large}, a.s. for all but finitely many $n$
there are at least $n^2$ times in the interval $[2^n,2^{n+1})$ for which
$|X_t|>(2^n)^{2/3} n^{-(2+\ep)}$.

From the formulas for $\phi_{n,m}$ and $Z_m$ it is easily seen (see also
Proposition~\ref{prop:FMT_size}) that for large $\gamma$, uniformly for all
$k$ 
\[
\P(Y_t > \gamma k^2 | X_t = -k) > c_1 \gamma^{-3/2} ,
\]
and so the probability that this fails for the above times when $|X_t|$ is
large is at most
\[
(1-c_1 \gamma^{-3/2})^{n^2} < \pexp{-c_1 n^2 \gamma^{-3/2}} .
\]
For $\gamma=n^{4/3-\ep}$ this is finitely summable and so a.s.\ for all but
finitely many $n$ there is some $2^n\le t_n < 2^{n+1}$ with
\[
|X_{t_n}| > (2^n)^{2/3} n^{-(2+\ep)} ,
\]
and
\[
Y_{t_n} > n^{4/3-\ep} X_{t_n}^2,
\]
and so
\[
\lim Y_{t_n} \frac{n^{8/3+3\ep}}{(2^n)^{4/3}} = \infty .
\]
\end{proof}

\begin{proof}[Proof of Theorem~\ref{thm:hull_growth}]
Since $|\barB_r| = \sum_{t<T_r} Y_t$, the first part follows from
Theorem~\ref{thm:layertime} and Lemma~\ref{lem:Y_small}, and the second
part follows from Theorem~\ref{thm:layertime} and Lemma~\ref{lem:Y_large}.
\end{proof}

%%%%%%%%%%%%%%%%%%%%%%%%%%%%
\section{Ball Volume Growth}  \label{sec:ball_vol}
%%%%%%%%%%%%%%%%%%%%%%%%%%%%

In order to get a lower bound on the $|B_r|$ we need to find (with good
probability) a large number of vertices within a short distance of the
root. When investigating the hull $\barB_r$ we found a number of steps at
which a free triangulation of an $m$-gon was added to the triangulation. We
therefore wish to estimate not only the number of vertices in a free
triangulation of a disc, but the number of such vertices that are close to 
the boundary.

\begin{defn}\label{def:height}
For triangulation $T$ of a disc and a vertex $v \in T$, the {\em height} of
$v$, denoted $h_v$, is the distance from $v$ to the boundary.
\end{defn}

Typically we expect a free triangulation of an $m$-gon to have size roughly
$m^2$. We will see that typically, most vertices in such a triangulation
have height at most roughly $\sqrt{m}$. This implies that typical distances
in $T$ are on the scale of $|T|^{1/4}$, conforming with the result of
\cite{ChSc}. The following methods may also be used to estimate the maximal
height in a free triangulation.

\begin{lemma}\label{lem:big_layer}
Let $T$ be a free triangulation of an $m$-gon. For any $\ep>0$ there are
$c_1,c_2>0$, such that for sufficiently large $m$
\[
\mu_m(A) > c_2 ,
\]
where $A$ is the event consisting of all triangulations, $T$, of an
$(m+2)$-gon with
\[
|\{ u \in T, h_u < \sqrt{m} \log^{3+\ep} m\}| > c_1 m^2 .
\]
\end{lemma}

Based on this Lemma, the proof of Theorem~\ref{thm:ball_growth} is
straightforward --- simply find a free triangulation of a large cycle near
the root and it will follow that there are many vertices near the root.

\begin{proof}[Proof of Theorem~\ref{thm:ball_growth}]
Fix $\ep>0$ and let
\[
L_r =
\big\{ t<T_r {\textrm\ s.t.\ } |X_t| > r^2 \log^{-(6+5\ep/3)} r \big\} .
\]

Using the lower bound on $T_r$ (Theorem~\ref{thm:layertime}) together with
Lemma~\ref{lem:X_large} to find times when $|X_t|$ is large, we get that
for some $c_1$, a.s.\ for all but finitely many $r$
\[
|L_r| > c_1\log^2 r.
\]

At each time $t\in L_r$ a free triangulation is glued to a cycle, with
every vertex in the cycle at distance at most $r$ from the root. For each
such $t$ the probability that the event in Lemma~\ref{lem:big_layer} fails
to occur for this free triangulation is at most $c_2$. Thus the probability
of failure at all $t\in L_r$ is at most $\pexp{-c_2\log^2 r}$ for some
$c_2>0$, which is finitely summable. Thus a.s.\ for all but finitely many
$r$ the event of Lemma~\ref{lem:big_layer} occurs at some time $t\in L_r$.

By the definition of $L_r$ and from the Lemma, this implies that there are
at least $r^4 \log^{-(12+10\ep/3)} r$ vertices at bounded distance from the
cycle, and thus also at bounded distance from the root. To find this bound,
note that $|X_t|<M_t$ and $t<T_r$. Theorem~\ref{thm:M_small} bounds $M_t$
in terms of $t$ and Theorem~\ref{thm:layertime} bounds $T_r$ in terms of
$r$. Combined they give
\[
|X_t| < c_3 r^2 \log^3 r ,
\]
hence the bound on height in Lemma~\ref{lem:big_layer} is
\[
\sqrt{|X_t|} \log^{3+\ep} |X_t| < c_4 r \log^{9/2+\ep} r .
\]

Summarizing, for all but finitely many $r$ we found
$r^4 \log^{-(12+10\ep/3)} r$ vertices in $B_{r'}$ with
$r' = r+c_4 r \log^{9/2+\ep} r$. By monotonicity of $|B_r|$ this suffices.
\end{proof}

In proving Lemma~\ref{lem:big_layer} we make use of the following distribution,
closely related to the free triangulation of a disc.

\begin{defn}
The free {\em marked} triangulation of a disc is a distribution on
triangulations of the disc with a marked internal vertex, that assigns a
rooted triangulation $T$ marked at $v$ probability. For an $(m+2)$-gon:
\[
\tilde \mu_m(T,v) = \tilde Z_m \alpha^{-|T|},
\]
where
\[
\tilde Z_m = \sum n \phi_{n,m} \alpha^{-n}.
\]
\end{defn}

This may also be defined as the annealed distribution for marking a random
vertex in a free triangulation of the disc.
Proposition~\ref{prop:free_size} implies 
\[
\tilde Z_m = Z_m \E_{\mu_m}|T| =
\frac16 \binom{2m+2}{m} \left(\frac94\right)^{m+1} ,
\]
hence the relation between $\tilde \mu_m$ and $\mu_m$ may be written as
\[
\tilde \mu_m(T) = \frac{3|T|}{(m+1)(2m+1)} \mu_m(T) .
\]
Clearly, conditioned on the size of $T$, the triangulation marginal of
$\tilde \mu_m$ and $\mu_m$ are equal, and conditioned on the triangulation, 
each internal vertex has probability $|T|^{-1}$ of being the marked one.

The size distribution of the free marked triangulation is interesting, as
is evident from the following proposition.

\begin{prop} \label{prop:FMT_size}
The distributions of $m^{-2}|T|$ with respect to $\tilde \mu_m$ converge to
that of $2/3S_{1/2}$, where $S_{1/2}$ is an asymmetric stable random
variable.
\end{prop}

$S_{1/2}$ with the Levi distribution, is one of very few cases of stable
random variables for which there is a (more or less) closed form for the
density function. $S_{1/2}$ has the same law as $g^{-2}$ where $g$ is a
standard Gaussian random variable. In particular, for large $m$, the size
of a free marked triangulation of an $m$-gon is distributed approximately
as $2/3 g^{-2}$. As a consequence, the size of a free (unmarked)
triangulation of an $m$-gon also converge in distribution.

\begin{proof}
The probability of a free marked triangulation having size $n$, namely
\[
\tilde \mu_m(|T|=n) = \frac{n\phi_{n,m}}{\tilde Z_m \alpha^{-n}}
\]
may be rewritten as 
\[
\frac{8n(m+2)}{3(2n+2m+1)(2n+2m+2)} \cdot
\left(\left(\frac4{27}\right)^n \binom{3n}{n}\right) \cdot
\prod_{i=1}^{2m} \frac{n+i/3}{n+i/2} .
\]
For large $n$, the second term is roughly $\sqrt{\frac{3n}{4\pi}}$. For
$n \gg m \gg 1$, the first term is roughly $2m/3n$ while the last is
roughly $e^{-m^2/3n}$. These approximations are uniform for all
$n \gg m \gg 1$. It follows that
\[
\tilde \mu_m(|T|=n) \sim \frac{mn^{-3/2}}{\sqrt{3\pi}}e^{-m^2/3n} .
\]
And by a change of variable, for any $t$
\[
\lim_{m \to \infty} \tilde \mu_m(|T| > tm^2) =
  \frac{1}{\sqrt{3\pi}}\int_t^\infty x^{-3/2} e^{-1/3x} dx ,
\]

Hence the distribution of $m^{-2}|T|$ with respect to $\tilde \mu$
converges as $m \to \infty$ to that of $2/3S_{1/2}$.
\end{proof}

\begin{figure}
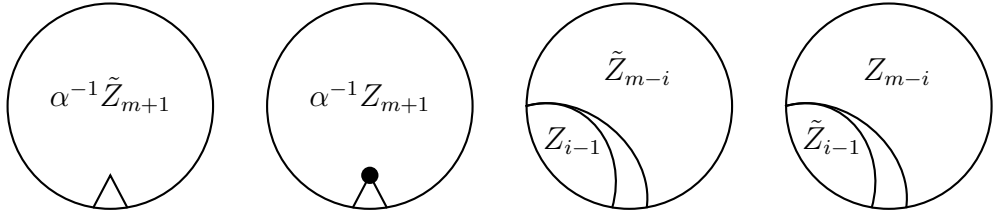

\begin{center}
\psset{unit=1.15mm}
\pspicture(-12,-12)(102,12)
  \pscircle(0,0){12}
  \psline(-2,-11.8)(0,-8)(2,-11.8)
  \put(-7,0){$\alpha^{-1}\tilde Z_{m+1}$}

  \pscircle(30,0){12}
  \psline(28,-11.8)(30,-8)(32,-11.8)
  \pscircle*(30,-8){1}
  \put(23,0){$\alpha^{-1}Z_{m+1}$}

  \pscircle(60,0){12}
  \psbezier(58,-11.8)(60,-5)(55,2)(48,0)
  \psbezier(62,-11.8)(63,-5)(55,2)(48,0)
  \put(57,3){$\tilde Z_{m-i}$}
  \put(50,-5){$Z_{i-1}$}

  \pscircle(90,0){12}
  \psbezier(88,-11.8)(90,-5)(85,2)(78,0)
  \psbezier(92,-11.8)(93,-5)(85,2)(78,0)
  \put(87,3){$Z_{m-i}$}
  \put(80,-5){$\tilde Z_{i-1}$}
\endpspicture
\end{center}
\caption{\label{fig:marked_peel}
The possibilities when a triangle is added in the free marked triangulation.}
\end{figure}

The free marked triangulation has the advantage over the free triangulation
that it may be sampled via a peeling process similar to the peeling process
for the UIPT. As each triangle is added, there are three possibilities,
shown in Figure~\ref{fig:marked_peel}. The triangle may partition the disc
in two, in which case there is a marked and unmarked triangulations in the
two parts, corresponding to the infinite and free parts in the construction
of the UIPT. The second option is that a new vertex is added. The third is
that the added vertex is the marked vertex.

As before, in order to analyze the sampling process, assume that triangles
are added to the boundary at in order, adding all triangles incident to a
vertex before moving to the next vertex counterclockwise. As before, this
does not effect the resulting measure or the distribution of the number of
required steps, but will make understanding the height distribution easier.

As before, we consider the size of the boundary, $M_n$, after $n$ steps
were taken. $M_0=m$ and at each step we have $M_{n+1}=M_n+X_n$ where the
distribution of $X_n$ is as in Equations~(\ref{eq:infnew}),
(\ref{eq:infold}) with $\tilde Z_m$ 
taking the role of $C_m$:
\bi
\item $X_t=1$ and the new vertex is unmarked with probability
$  \frac{\tilde Z_{m+1}}{\alpha \tilde Z_m}.  $
\item $X_t=1$ and the new vertex is marked with probability
$  \frac{Z_{m+1}}{\alpha \tilde Z_m}.  $
\item $X_t=-k$ with probability
$  \frac{\tilde Z_{m-k}Z_{k-1}}{\tilde Z_m}.  $
\ei

It is easy to see that as $M \to \infty$ these converge to the same
probabilities as for the UIPT. However, here the drift is reversed: the
small $M_n$'s are likelier to give negative $X$. Here the expectation of
the steps is deceiving, as it is positive. Still, the probabilities for
large negative steps are sufficient to cause the Markov chain to terminate
within roughly $M_0^{3/2}$ steps. Since each layer contains roughly $M_0$
vertices, this implies that the height of the marked vertex is at most on
the order of $\sqrt {M_0}$.

\begin{lemma}\label{lem:M2_small}
For $M_n$ as above, for some $c>0$
\[
\P(\max M_n > \lambda M_0) < c \lambda^{-1}.
\]
\end{lemma}

\begin{proof}
By Proposition~\ref{prop:FMT_size}, for a free marked triangulation $T$ of
an $(m+2)$-gon, 
\[
\tilde \mu_m(|T| > tm^2) < c_1t^{-1/2},
\]
and that for some $c_2>0$
\[
\tilde \mu_m(|T| > m^2) > c_2.
\]

Let $M_n$ be the sequence observed in the process of sampling of a
triangulation $T$ of an $(M_0+2)$-gon. Let $A$ be the event that for some
$n$, $M_n > \lambda M_0$. Conditioned on $A$, $T$ includes a free marked
triangulation of an $(M_n+2)$-gon, and so 
\[
\P\left(|T| > \lambda^2 M_0^2 \mid A\right) > c_2.
\]
Thus
\[
\P(A) < \frac{\P(|T| > \lambda^2 M_0^2)}{c_2} < \frac{c_1}{c_2 \lambda}.
\]
\end{proof}

\begin{lemma} \label{lem:s_small}
For $M_n$ as above, let $s$ be the smallest $n$ such that either
$M_n \le M_0/2$ or the process terminates at time $n$. For any $m\ge M_0$
and $\ep>0$ there is a $c=c(\ep)$ such that
\[
\P(s > M_0^{3/2} \log^{2+3\ep} m) < c\log^{-(1+\ep)} m.
\]
\end{lemma}

As before, let $\F_n$ denote the sigma field generated by the random
variables $M_0,\ldots,M_n$. Expectation with respect to $\F_n$ is denoted
by $\E_{\F_n}$ (expectation conditioned on the past).

\begin{proof}
This is similar to the proof of Theorem~\ref{thm:M_large}.
To avoid difficulties arising from the termination of the process, continue
the Markov chain after termination with $X_t=-1$ after 0 is hit. With this
continuation the following bounds hold universally. If we find that
$M_n<M_0/2$, then clearly $s<n$ whether or not the process terminated.

We show that at each step there is a not too small probability of the
Markov chain taking a step into $[0,M_0/2)$, and consequently this happens
after a small number of steps. Direct inspection of the step probabilities
reveals that for some universal $c_1$: 
\[
\P(X_n=-k \mid M_n) > c_1 (M_n-k)^{-1/2} M_n^{1/2} k^{-5/2}.
\]
Since $k<M_n$, it follows that the probability of taking a step from $M_n$
to the interval $[0,M_0/2)$ satisfies
\[
\P(M_{n+1}<M_0/2 \mid M_n) > c_2 M_0^{1/2} M_n^{-2}.
\]

Define the stopping time $t$ as the first time the Markov chain leaves the
interval $[M_0/2,M_0\log^{1+\ep} m]$ (or terminates). As long as $M_n$ is
in the interval the probability of leaving the interval by having
$M_{n+1}<M_0/2$ is at least $c_2 M_0^{-3/2} \log^{-(2+2\ep)} m$. It follows
that
\begin{eqnarray*}
\P(t>M_0^{3/2} \log^{2+3\ep} m)
  &<& \left( 1-c_2 M_0^{-3/2}\log^{-(2+2\ep)} m
      \right)^{M_0^{3/2} \log^{2+3\ep} m}               \\
  &<& \pexp{-c_2 \log^{\ep} m}.
\end{eqnarray*}

Clearly $s=t$ unless the process left the interval with
$M_t>M_0\log^{1+\ep} m$. By Lemma~\ref{lem:M2_small} the probability of
that is bounded by $c_3\log^{-(1+\ep)} m$, and therefore
\[
\P(s>M_0^{3/2} \log^{2+3\ep} m) <
   c_3\log^{-(1+\ep)} m + \pexp{-c_2 \log^{\ep} m}
\]
as claimed.
\end{proof}

\begin{lemma} \label{lem:hv_small}
Let $(T,v)$ be a free marked triangulation of an $(m+2)$-gon. For any
$\ep>0$, for all large enough $m$ 
\[
\tilde \mu_m(h_v > \sqrt m \log^{3+\ep} m) < 1/2.
\]
\end{lemma}

\begin{proof}
Consider for $k<\log_2 m$ the times $S_k$ defined as the smallest $s$ such
that $M_s \le 2^{-k}m$. We show that $S_{k+1}-S_k$ is unlikely to be large
and give a bound on the distance between the boundary at time $S_{k+1}$ and
the boundary at time $S_k$.

Using the Lemma~\ref{lem:s_small} applied to the process starting at time
$S_k$ with boundary size $M_{S_k} \le 2^{-k}m$, and with $m$ in the lemma's
formulation being the boundary size of the original disc, we get the bound
\[
\P(S_{k+1}-S_k > (2^{-k}M_0)^{3/2} \log^{2+\ep} m) < c\log^{-(1+\ep/3)} m.
\]
Therefore with probability at least $1-c\log^{\ep/3} m$, for all relevant
$k<\log_2 m$
\begin{equation}\label{eq:fast_half}
S_{k+1}-S_k < (2^{-k}m)^{3/2} \log^{2+\ep} m .
\end{equation}
For large $m$ this probability is at least 3/4.

As before, let $T_r$ denote the time when the $r$'th layer is complete,
i.e.\ the smallest $t$ such that the distance from the original boundary to
the boundary at time $t$ is $r$. At time $T_{r+1}$ all vertices on the
boundary have height $r+1$ and so they were all added after time $T_r$.
Since vertices are added at the boundary one at a time, we have
$T_{r+1}-T_r \geq M_{T_{r+1}}$.

If (\ref{eq:fast_half}) holds, then the number of rounds completed between
$S_k$ and $S_{k+1}$ is at most 
\[
1+\frac{S_{k+1}-S_k}{2^{-(k+1)}m} < 1+2\sqrt{m} \log^{2+\ep} m,
\]
and by summing up, with probability at least 3/4, the total number of
rounds completed before the process either terminates or reaches $M_n=0$ is
at most $c_1 \sqrt{m} \log^{3+\ep} m$. 

If the process terminated before $M_n=0$, then
$h_v<c_1 \sqrt{m} \log^{3+\ep} m$ as required. In the case that the process
did not terminate and reached $M_n=0$, in the remaining 2-gon there is a
free marked triangulation, and with probability 3/4 the number of vertices
inside is at most 30 and so $\P(h_v>30+c_1\sqrt{m} \log^{3+\ep} m) < 1/2$. 
By changing $\ep$ the constants $c_1$ and 30 may be disposed of.
\end{proof}

\begin{proof}[Proof of Lemma~\ref{lem:big_layer}]
Let $(T,v)$ be a free marked triangulation of an $m$-gon. Define $t=|T|$
and
\[
s = |\{ u \in T | h_u < \sqrt m \log^{3+\ep} m\}|.
\]
The probability that $h_v<\sqrt m \log^{3+\ep} m$, bounded by
Lemma~\ref{lem:hv_small}, is clearly the expectation of $s/t$ with respect
to $\tilde \mu_m$, i.e.
\[
\sum_T \frac{s}{t} \tilde\mu_m(T) > 1/2.
\]
Using the relation between $\tilde \mu_m$ and $\mu_m$, this translates to
\[
\sum_T \frac{s}{t} \frac{3t}{(m+1)(2m+1)} \mu_m(T) > 1/2,
\]
or
\[
\E_{\mu_m} s = \sum_T s \mu_m(T) > \frac{m^2}{3}.
\]

To finish, note that since $s<t$ and
$\mu_m(t > \gamma m^2) < c_1\gamma^{-3/2}$, the expectation of $s$
restricted to any event of probability $p$ is bounded by $c_2m^2p^{1/3}$
for some $c_2$. Therefore, if the probability that $s<c_3 m^2$ is $p$, the
expectation of $s$ is at most 
\[
pc_3 m^2 + c_2m^2(1-p)^{1/3}.
\]
For any $c_3<1/3$ this gives a lower bound on $p$ independent of $m$.
\end{proof}

%%%%%%%%%%%%%%%%%%%%%%%%%%%%%%%%%
\section{Percolation on the UIPT}  \label{sec:percolation}
%%%%%%%%%%%%%%%%%%%%%%%%%%%%%%%%%

The peeling construction may also be used to understand percolation on the 
UIPT. Bernoulli site percolation with parameter $p$ is the probability
measure on colorings of the graph's vertices where each vertex is colored
black with probability $p$ and white otherwise, independently of all other
vertices. Percolation is defined as the event that the root vertex is in an
infinite connected black component. The resulting measure is denoted by
$\P_p$.

Since the graph in question, namely the UIPT, is a random graph there is a
distinction between annealed and quenched statements about percolation
(though since there is no interaction between the percolation and the
underling graph the underlying probability measure is the same). Annealed
statements are averaged on all planar triangulations. Quenched statements
on percolation are about the coloring conditioned on the triangulation.
In particular, the critical probability $p_c$ for annealed percolation is
the infimum of all $p$ such that $\P_p(\textrm{percolation})>0$. This means
that in a positive measure of triangulations (w.r.t.\ $\tau$) the
probability of percolation is positive. The quenched problem is what is the
least $p$ such that on $\tau$-almost all triangulations there is a positive
probability of percolation. The answers turn out to be the same.

\begin{thm}\label{thm:annperc}
The annealed critical probability for site percolation on the UIPT is 1/2.
\end{thm}

\begin{proof}
Suppose each vertex is colored independently at random, black with
probability $p$ and white otherwise. Assume with out loss of generality
that the root vertex is black. We sample the connected black component
containing the root vertex by adding triangles one by one to the
triangulation, as in the peeling construction of the UIPT. As new vertices
are added we color them randomly, independently of all previous events.
This samples site percolation on (a part of) the UIPT.

If at any time all the vertices on the outer boundary of the sampled
sub-triangulation we see are white, then we have found a cycle of white
vertices enclosing the root, and so the black connected cluster containing
the root is necessarily finite. On the other hand, if the process continues
indefinitely and never reaches an all white outer boundary, then the root
is in an infinite connected black component. 

Recall that we are free to choose the edge in the outer boundary to which
we attach a triangle on at any time. At all times we will choose a boundary 
edge that has one black and one white endpoint (unless the outer boundary
of the sub-triangulation is monochromatic). Choosing the edge in this
manner will guarantee two things. First, the black (resp. white) vertices
on the outer boundary lie on a continuous arc along the boundary as in
Figure~\ref{fig:perc}(a,b). Second, since each new vertex is connected to
existing vertices of both colors (as long as the boundary is not
monochromatic), all black (resp. white) vertices form a single connected
cluster.

\begin{figure}
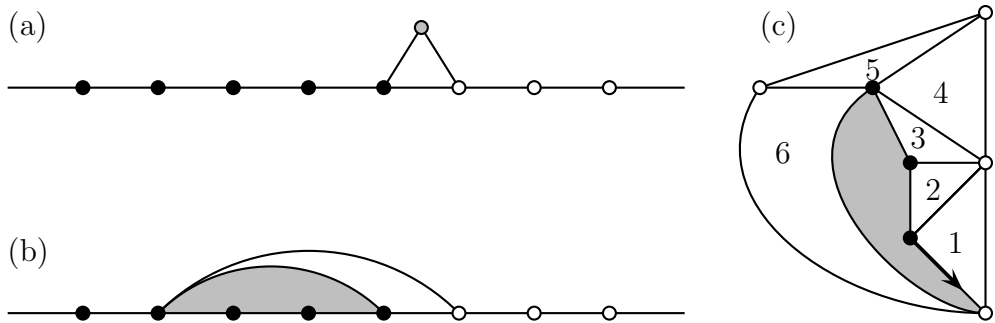

\begin{center}
\psset{unit=1mm}
\pspicture(0,0)(130,40)
  \put(0,7){(b)}
  \psarc[fillstyle=solid,fillcolor=lightgray](35,-15){21.213}{45}{135}
  \psarc(40,-20){28.284}{45}{135}
  \psline(0,0)(90,0)
  \multips(10,0)(10,0){5}{\pscircle*{1}}
  \multips(60,0)(10,0){3}{\pscircle[fillstyle=solid]{1}}

  \put(0,37){(a)}
  \psline(50,30)(55,38)(60,30)
  \pscircle[fillstyle=solid,fillcolor=lightgray](55,38){1}
  \psline(0,30)(90,30)
  \multips(10,30)(10,00){5}{\pscircle*{1}}
  \multips(60,30)(10,00){3}{\pscircle[fillstyle=solid]{1}}

  \put(100,37){(c)}
  \psbezier(130,0)(110,0)(90,15)(100,30)
  \psbezier[fillstyle=solid,fillcolor=lightgray]
           (130,0)(120,0)(100,20)(115,30)
  \psline[fillstyle=solid](130,40)(100,30)(115,30)(120,20)(120,10)
         (130,0)(130,40)(115,30)(130,20)(120,10)(130,20)(120,20)
  \psline[linewidth=.5]{->}(120,10)(127,3)
  \pscircle[fillstyle=solid,fillcolor=white](130,40){1}
  \pscircle[fillstyle=solid,fillcolor=white](100,30){1}
  \pscircle*(115,30){1}
  \pscircle*(120,20){1}
  \pscircle[fillstyle=solid,fillcolor=white](130,20){1}
  \pscircle*(120,10){1}
  \pscircle[fillstyle=solid,fillcolor=white](130,0){1}
  \put(125,8){1}
  \put(122,15){2}
  \put(120,22){3}
  \put(123,28){4}
  \put(114,31){5}
  \put(102,20){6}
\endpspicture
\end{center}
\caption{\label{fig:perc} Sampling annealed percolation on the UIPT.}
\end{figure}

When a new vertex is added, we color it randomly (as in
Figure~\ref{fig:perc}(a), showing a segment of the boundary). When the new
triangle includes only old vertices, as in Figure~\ref{fig:perc}(b), the
event of percolation does not depend on the finite part of the
triangulation (in the shaded region). Thus, to determine whether
percolation occurs we only need to keep track of the outer boundary and of
the number of black and white vertices in it. 

Figure~\ref{fig:perc}(c) shows a possible outcome of the first few steps.
The triangles are numbered in the order they are added. After the sixth
triangle is added, we know that the black cluster of the root is finite.
Since the part of the triangulation enclosed in the hole is finite,
percolation can not occur.

\medskip

Let $B_n$ (resp. $W_n$) denote the number of black (resp. white) vertices
on the boundary at time $n$. As before, let $M_n+2=B_n+W_n$ denote the size
of the boundary after $n$ steps. Also let $X_n$ denote the step size, with
distribution given by (\ref{eq:step_dist}). If $X_n=1$, then the new vertex
is colored black or white with the given bias. Therefore conditioned on
$(B_n,W_n)$:
\[
(B_{n+1},W_{n+1}) = \left\{ \begin{array}{ll}
(B_n+1,W_n\phantom{+1})  &  \textrm{with prob. $p$}         \\
(B_n\phantom{+1},W_n+1)  &  \textrm{with prob. $1-p$}       \\
\end{array} \right. \quad \textrm{If $X_n=1$.}
\]

On the other hand, if $X_n<0$, then $|X_n|$ vertices from the outer boundary
become internal vertices. In this case it is equally likely that the
removed vertices are on the left or the right of the new triangle (the
triangulation's geometry is independent of previously chosen colors). Note
however, that if $B_n+X_n<0$ and the removed vertices are on the black
side, this does not result in a negative number of black vertices on the
boundary. Instead, all the black vertices, and some ``borrowed'' white ones
will be removed. Therefore conditioned on $(B_n,W_n)$:
\[
(B_{n+1},W_{n+1}) = \left\{ \begin{array}{ll}
f(B_n+X_n,W_n\phantom{+X_n})  &  \textrm{with prob. $1/2$}      \\
f(B_n\phantom{+X_n},W_n+X_n)  &  \textrm{with prob. $1/2$}      \\
\end{array} \right.\quad \textrm{If $X_n<0$,}
\]
where $f(B,W)=(B,W)$ unless $B<0$ or $W<0$:
\[
f(B,W) = \left\{ \begin{array}{ll}
(B,W)   & \textrm{if $B,W\ge0$,}        \\
(0,B+W) & \textrm{if $B<0$,}        \\
(B+W,0) & \textrm{if $W<0$.}        \\
\end{array} \right.
\]

We see that if $p=1/2$, then $X_n$ is added with equal probabilities to $B$
or $W$, but if $p\neq 1/2$ there is a bias, and since $\E X_n \to 0$, the
less probable color tends to die out, as we now prove.

\medskip

Suppose $p<1/2$. Consider the Markov chain $(B'_n,W'_n)$ not limited to
positive values, with evolution
\[
(B'_{n+1},W'_{n+1}) = \left\{ \begin{array}{ll}
(B'_n+1,W'_n\phantom{+1})  &  \textrm{with prob. $p$}           \\
(B'_n\phantom{+1},W'_n+1)  &  \textrm{with prob. $1-p$}     \\
\end{array} \right. \quad \textrm{If $X_n=1$.}
\]
\[
(B'_{n+1},W'_{n+1}) = \left\{ \begin{array}{ll}
(B'_n+X_n,W'_n\phantom{+X_n})  &  \textrm{with prob. $1/2$}     \\
(B'_n\phantom{+X_n},W'_n+X_n)  &  \textrm{with prob. $1/2$}     \\
\end{array} \right.\quad \textrm{If $X_n<0$,}
\]
coupled in the natural way to $(B_n,W_n)$. The difference $B'_{n+1}-B'_n$
is equal to the $B_{n+1}-B_n$ except in two cases:
\bi
\item $W_n+X_n<0$ and additional black vertices are removed. In this event
$B'_{n+1}-B'_n > B_{n+1}-B_n$.
\item $B_n+X_n<0$ and additional white vertices are removed. In this event
percolation does not occur.
\ei
Therefore, conditioned on the event of percolation, $B_n \le B'_n$. We will 
see that a.s.\ for some $n$, $B'_n \le 0$ and so
$\P_p(\textrm{percolation})=0$. Denote by $\F_n$ the $\sigma$-field
generated by the random variables $B_0,W_0,\ldots,B_n,W_n$. Note that $M_n$
and $B'_n$ are determined by $B,W$ up to time $n$, and hence are
$\F_n$-measurable.
\[
\E_{\F_n} B'_{n+1}-B'_n \leq 1/2\E_{\F_n} X_n - (1/2-p) \P(X_n=1 | \F_n) .
\]

Recall that if $m>m'$, then the distribution of $X_n$ conditioned on $M_n=m$ 
is stochastically dominated by the distribution conditioned on $M_n=m'$.
Since a.s.\ $M_n \to \infty$, from some time on $M_n$ is larger than any
fixed $M$. Since $\E_{\F_n} X_n \sim cM_n^{-1/2} \to 0$ from some time on
$B'_n$ is dominated by a random walk with steps bounded from above
($X_n\le1$) and negative expectation. Such a random walk will a.s.\ tend to
$-\infty$, and since conditioned on percolation $B'_n$ remains positive,
$\P_p(\textrm{percolation}) = 0$. 

Moreover, if whenever $B_n=0$ we ``reset'' by setting $B'_n=0$, then
$B_n \le B'_n$ holds unconditionally, and when $B'_n$ is large enough 
it still has negative expected change. With this modification $B'_n$ does
not tend to $-\infty$ but has a stationary distribution with exponential
decay. By Borel-Cantelli it follows that $B_n \le B'_n = O(\log n)$.

When $p>1/2$, the roles of $B_n$ and $W_n$ are interchanged, and since
$B_n+W_n$ grows like $n^{2/3}$ it follows that a.s.\ $B_n \to \infty$. In
particular with positive probability at all times $B_n>0$, i.e.,
$\P_p(\textrm{percolation}) > 0$.  
\end{proof}

Thus, we see that if $p<1/2$, then a.s.\ the black clusters are finite,
while if $p>1/2$, then with positive probability there is an infinite black
cluster (since white clusters are all finite, it is unique). To see that
this is a.s.\ the case we prove the following 0-1 law. A general 0-1 law
for triangulations is proved in \cite{BCR}. However, we also need a 0-1 law
for annealed percolation on the UIPT.

\begin{thm} \label{thm:zero_one}
The probability of any event invariant to finite changes in the
triangulation is 0 or 1. Moreover, the same is true for the annealed
percolation on the UIPT with any $p$.
\end{thm}

\begin{proof}[Proof of Theorem~\ref{thm:percolation}]
$p_c=1/2$ is an event independent of finite changes in the UIPT. Since the
annealed $p_c=1/2$, $\tau(p_c=1/2)>0$.

For $p=1/2$ we see that the probability of having an infinite black (or, by
symmetry, white) cluster is either 0 or 1. Thus, we see that either a.s.\
there are no monochromatic infinite clusters, or a.s.\ there are infinite
clusters of both colors. To rule out the latter possibility note that there
are a.s.\ infinitely many times when $|X_t|>M_t/2$. Conditioned on
$|X_t|>M_t/2$, with probability at least 1/2 one of the colors is
completely removed from the outer boundary. Thus, a.s.\ one of the colors
dies out infinitely often and has only finite clusters. By symmetry,
$\P_{1/2}(\textrm{percolation})\le 1/2$, and hence it is 0.  
\end{proof}

\begin{proof}[Proof of Theorem~\ref{thm:zero_one}]
We represent the UIPT as a function of an infinite sequence of independent
random variables, such that a.s.\ changing a finite number of them will
only change a finite sub-triangulation. Since any event determined by the
tail of such a sequence has probability 0 or 1, this is sufficient.

The basis will be the peeling construction of the UIPT. For each $m,n,i$ 
define an independent random variable $Z_{m,n,i}$, with the appropriate
distribution for a step when the boundary size is $m$. Such a random
variable includes which is the third vertex of an added triangle as well as
a sample of the free triangulation when appropriate. $Z_{m,n,i}$ is used at
time $t$ under 3 conditions:
\bi
\item $M_t=m$
\item $\max_{s<t} M_s = n$
\item If $s$ is the first time that $M_s=n$, then $t-s=i$.
\ei

For example, if the observed sequence of $M$'s is $1,2,3,1,2,3,1,2$, then
the next step uses $Z_{2,3,5}$: $M_t=2$, the maximal value seen so far is
3, first reached 5 steps ago.

Since a.s.\ $M_t \to \infty$, changing a finite number of the $Z_{m,n,i}$
will only change the evolution of the triangulation in a finite number of
steps. Once $M_t$ is sufficiently large the evolution will not be changed at
all.

\medskip

The exact same proof also works for percolation on the UIPT, adding the
random colors of new vertices to the variables $Z_{m,n,i}$.
\end{proof}

\bigskip

\filbreak
\begingroup
\small
\parindent=0pt

\vtop{
\hsize=2.3in
Omer Angel\\
Department of Mathematics\\
Weizmann Institute of science\\
Rehovot, 76100, Israel\\
{omer@wisdom.weizmann.ac.il}
}
\endgroup
\filbreak

\end{document}